\documentclass[10pt,b5paper]{article}
\usepackage[T1]{fontenc}
\usepackage[utf8]{inputenc}
\usepackage[english]{babel}
\usepackage{amsmath,amssymb,amsthm,cite}
\usepackage{lmodern}
\usepackage{geometry}
\usepackage{graphicx}
\usepackage{url}
\usepackage{titlesec}

\titleformat{\paragraph}[runin]{\normalfont\itshape}{}{0pt}{}[.]
\titlespacing*{\paragraph}{0pt}{0.7em}{0.5em}

\newtheorem{theorem}{Theorem}[section]
\newtheorem{lemma}[theorem]{Lemma}
\newtheorem{proposition}[theorem]{Proposition}
\newtheorem{corollary}[theorem]{Corollary}

\theoremstyle{definition}
\newtheorem{definition}[theorem]{Definition}
\theoremstyle{remark}
\newtheorem{remark}[theorem]{Remark}
\numberwithin{equation}{section}

\DeclareMathOperator{\E}{E}
\DeclareMathOperator{\Var}{Var}
\DeclareMathOperator{\Cov}{Cov}
\DeclareMathOperator{\Corr}{Corr}
\DeclareMathOperator{\Bin}{Bin}
\DeclareMathOperator{\Mult}{Mult}
\DeclareMathOperator{\sgn}{sgn}

\newcommand{\Z}{\mathbb{Z}}

\newcommand{\hadj}{h_{\mathrm{adj}}}

\title{Discrete total variation under the multinomial,\\Part 5:\\
cyclic total variation}
\author{
Nikola Bani\'c\\
Gideon Brothers, 10000 Zagreb, Croatia\\
\texttt{nbanic@gmail.com}
\and
Neven Elezovi\'c\\
Department of Applied Mathematics,\\
Faculty of Electrical Engineering and Computing,\\
University of Zagreb, 10000 Zagreb, Croatia\\
\texttt{neven.elezovic@fer.hr}
}
\date{\today}

\begin{document}
\maketitle

\begin{abstract}
We study the cyclic discrete total variation obtained by adding the wrap-around
edge to the usual open-chain statistic. Under the uniform multinomial law the
cyclic symmetry recombines the same per-edge, adjacent, and disjoint covariance
blocks used for the linear statistic, but with different multiplicities, and this
simplifies the variance: the bin-independent leading coefficient \(2A+4B\) carries
over, the constant-order correction collapses \emph{exactly} to \(\beta_c(n)=\beta_1 n\)
(no \(n^0\) or \(1/n\) tail; obtained directly here, and already reached in
Part~2 through the transfer identity), and a conditional second-order coefficient
\(\gamma_c(n)\) is given. The exact finite-\(N\) cyclic--linear transfer identity of
Part~2, read order by order, both delivers these and, as its main consequence for the
\emph{linear} statistic, explains the traveling mode of the linear variance: comparison
with the (first-order) mode-free symmetric cyclic profile shows that deleting one edge
introduces a rising contribution whose balance against the falling cyclic bulk creates a
ridge at \(n\asymp\sqrt{K_0N}\) with local optimizer \(K_0=C_1/\beta_1\). Uniform control
of three ridge blocks proves this optimizer on every compact ridge window, upgrading the
previously heuristic \(\sqrt N\) mode law to a theorem; the global localization of the
exact maximizer remains open. We also record the corresponding cyclic reference-centered
statistic under a fixed alternative and the leading Nabeya variance formula for the
matched non-uniform cyclic statistic. Throughout, ``exact'' refers to the
asymptotic coefficient in \(N\); the operative symmetry is cyclic stationarity
(rotation invariance), not exchangeability.
\end{abstract}



\section{Introduction}

	The four preceding papers \cite{banic_elezovic_p1,banic_elezovic_p2,banic_elezovic_p3,banic_elezovic_p4} study the linear discrete total variation
	$d=\sum_{i=1}^{n-1}|X_{i+1}-X_i|$, introduced as a histogram statistic
	in~\cite{banic_elezovic_tvor}, whose summation stops at the two ends of the
	histogram. Those two endpoints break an otherwise perfect symmetry: the interior edges
	are on an equal footing, but the first and last bins each touch only one neighbor. \emph{The
	present paper closes the histogram into a ring} and studies the \emph{cyclic} total
	variation
	\[
		d^c=\sum_{j=1}^{n}|X_{j+1}-X_j|,\qquad X_{n+1}:=X_1,
	\]
	in which the wrap-around edge restores full symmetry: all $n$ edges are now
	cyclically stationary: their joint law is invariant under the rotation group
	$\mathbb Z/n$ (not exchangeable in the strict sense --- the covariances depend on
	cyclic distance). Symmetry is not a cosmetic gain. It sharpens every result of the variance
	theory and, unexpectedly, it explains a feature of the \emph{linear} statistic that the
	linear analysis can only assume --- and it proves the comb test's mode constant $K_0$ (formed from the Part~2 constants $C_1,\beta_1$) on every compact ridge window, upgrading the previously heuristic $\sqrt N$ mode law to a theorem.
	In the notation used in the later papers of the series this is the cyclic
	reference statistic \(d_{\mathbf u}^{c}\) for the uniform law
	\(\mathbf u=(1/n,\dots,1/n)\), abbreviated here to \(d^c\).

	Under the uniform null the leading dispersion coefficient (the $N\to\infty$ limit) is \emph{exactly}
	independent of the number of bins; its constant-order coefficient is \emph{exactly}
	linear, $\beta_c(n)=\beta_1n$, with no $n^0$ or $1/n$ term (``exactly'' throughout
	referring to the coefficient's structure in $n$, not to finite-$N$ exactness;
	this collapse was already reached in Part~2 through the transfer identity, and is
	obtained directly here); the
	first-order dispersion \emph{profile} is monotone decreasing in $n$, so it has no
	interior maximum, and numerically the exact finite-$N$ dispersion likewise shows no
	interior mode. Proving this exact monotonicity requires a uniform remainder estimate and is left
	open;
	and the second-order coefficient $\gamma_c(n)$ has a conditional closed form, pending
	the same triangular-lattice constant that appears for the linear statistic. Each of
	these is cleaner than its linear counterpart, and the comparison isolates exactly what
	the endpoints cost.

	The link is an exact finite-$N$ identity from
	Part~2~\cite[Lemma~3.6, Eq.~(3.7)]{banic_elezovic_p2}: the linear statistic is
	obtained from the cyclic one by deleting one edge. Read order by order, this identity
	delivers each linear asymptotic constant from the cyclic one and the per-edge block;
	it also explains the \emph{mode} of the linear variance, the constant
	$K_0=C_1/\beta_1=16.0429\ldots$ of practical importance in the
	comb test~\cite{banic_elezovic_ct}. The symmetric cyclic profile has no mode at first order; deleting one edge introduces
	a rising contribution that, balanced against the falling cyclic bulk, creates a
	traveling ridge in the regime $n\asymp\sqrt N$ with the local optimizer $K_0$. We
	therefore develop here -- as a consequence of the cyclic transfer -- the asymptotic
	analysis of that ridge: the leading optimizer $K_0$, its sub-leading coefficient, and
	the uniform block estimates that make the local ridge expansion rigorous. The
	remaining global step is to prove that the exact integer maximizer cannot escape
	this ridge window. We also record the non-uniform cyclic statistic.

	These results form a series on the exact and asymptotic moment theory of discrete total
	variation under the multinomial model, organized by \emph{which} moment is computed and
	\emph{under which} null: the expectation and variance of the linear statistic under the
	uniform null and under a general distribution; and the cyclic statistic -- this paper.
	The exact multinomial covariances and correlations on which the analysis rests are
	recalled in the Preliminaries. Throughout the series, results not proved in the
	surrounding text are flagged as conditional, numerically confirmed, or conjectured; the
	cyclic chapter inherits two remaining unresolved points -- the global monotonicity of the
	exact finite-$N$ dispersion, and the triangular-lattice constant in $\gamma_c$.

\section{Preliminaries and cyclic structure}
	\label{sec:multinomial-basics}

	We first recall the multinomial identities needed later and then introduce the
	cyclic statistic. The purpose of this section is only to fix notation and to
	record the exact structural reductions used throughout the paper: edge
	covariances, cyclic multiplicities, and the transfer identity relating the
	cyclic and open-chain variances.

\subsection{Multinomial distribution: definition and basic properties}
	\label{ssec:mult-def}

	Let $X = (X_1, X_2, \ldots, X_n)$ be a random vector following the
	multinomial distribution
	\begin{equation}\label{eq:mult}
		X \sim \Mult(N;\, p_1, p_2, \ldots, p_n), \qquad
		\sum_{i=1}^n p_i = 1,\quad p_i \ge 0.
	\end{equation}
	We interpret the components $X_i$ as the number of observations falling into the $i$-th bin
	in $N$ independent trials with probabilities $p_i$. The associated
	probability mass function is
	\begin{equation}\label{eq:mult-pmf}
		\Pr(X_1 = k_1, \ldots, X_n = k_n) =
		\frac{N!}{k_1!\,k_2!\cdots k_n!}\,p_1^{k_1}p_2^{k_2}\cdots p_n^{k_n},
		\quad \sum_i k_i = N,
	\end{equation}
	and $0$ otherwise.

	\paragraph{Marginal distribution of a single bin} The marginal distribution of a
	single bin is binomial,
	\begin{equation}\label{eq:mult-marg-1}
		X_i \sim \Bin(N, p_i),
	\end{equation}
	so $\E[X_i] = N p_i$ and $\Var(X_i) = N p_i(1 - p_i)$.

	\paragraph{Covariance structure} The standard computation (via the representation
	$X_i = \sum_{k=1}^N \mathbf 1[Z_k = i]$ with $Z_k$ iid categorical on
	$\{1, \ldots, n\}$ with probabilities $p_i$) gives
	\begin{equation}\label{eq:mult-cov}
		\Cov(X_i, X_j) =
	\begin{cases}
		N p_i (1 - p_i), & i = j, \\
		- N p_i p_j, & i \ne j.
	\end{cases}
	\end{equation}

\subsection{Edge variables and their basic covariances}
	\label{ssec:edge-vars}

	The total variation is built from the differences between neighboring bins, so
	these differences -- which we call \emph{edge variables} -- are the true objects
	of study; the individual counts $X_i$ will scarcely appear again on their own.
	Before any absolute value is taken, their joint behavior is completely
	elementary, and we record it explicitly. The lemma below is the
	last purely classical result we will need: from here on, the absolute values
	enter and the analysis becomes specific to the present problem.

\begin{definition}[Edge variables]\label{def:Dj}
	For a multinomial vector $X = (X_1, \ldots, X_n)$ we introduce the \emph{edge variables}
	\begin{equation}\label{eq:Dj-def}
		D_j := X_{j+1} - X_j, \qquad j = 1, \ldots, n-1,
	\end{equation}
	and their centered versions
	\begin{equation}\label{eq:tilde-Dj-def}
		\tilde D_j := D_j - \E[D_j] = D_j - N(p_{j+1} - p_j).
	\end{equation}
	Under the uniform law ($p_i \equiv 1/n$), $\E[D_j] = 0$ so $\tilde D_j = D_j$.
	In the non-uniform case, $\E[D_j] \ne 0$ in general.
\end{definition}

	The covariance structure of the edge variables (\emph{without absolute values})
	is exact and completely determined by (\ref{eq:mult-cov}) and the bilinearity
	of the covariance.

\begin{lemma}[Exact edge covariances]\label{lem:cov-Dj-multinomial}
	Under the multinomial $\Mult(N; p_1, \ldots, p_n)$:
	\begin{align}
		\Var(D_j) &= N\bigl[(p_j + p_{j+1}) - (p_j - p_{j+1})^2\bigr], \label{eq:Var-Dj} \\
		\Cov(D_j, D_{j+1})
		&= N\bigl[-p_{j+1}(1 + p_j + p_{j+2} - p_{j+1})
		+ p_j p_{j+2}\bigr], \label{eq:Cov-Dj-adj} \\
		\Cov(D_j, D_k)
		&= -N(p_{j+1} - p_j)(p_{k+1} - p_k),
		\qquad |k - j| \ge 2. \label{eq:Cov-Dj-disj}
	\end{align}
	Centering changes neither the variance nor the covariance: $\Var(\tilde D_j) = \Var(D_j)$,
	$\Cov(\tilde D_j, \tilde D_k) = \Cov(D_j, D_k)$.
\end{lemma}
\noindent
	The adjacent form (\ref{eq:Cov-Dj-adj}) appears asymmetric in $j$ and $j+2$
	only because $X_{j+1}$ is the shared variable.

\begin{proof}
	All three claims follow from the bilinearity of the covariance: for random
	variables $A, B, C, D$
	\begin{equation}\label{eq:cov-bilin}
		\Cov(A - B,\, C - D) = \Cov(A, C) - \Cov(A, D) - \Cov(B, C) + \Cov(B, D),
	\end{equation}
	and from the multinomial covariances (\ref{eq:mult-cov}).

	\medskip
	\emph{Proof for $\Var(D_j)$.} Special case of (\ref{eq:cov-bilin}) with
	$A = C = X_{j+1}$, $B = D = X_j$:
	\[
		\Var(D_j) = \Var(X_{j+1}) + \Var(X_j) - 2\Cov(X_j, X_{j+1}).
	\]
	Substituting (\ref{eq:mult-cov}):
	\begin{align*}
		\Var(D_j)
		&= N p_{j+1}(1-p_{j+1}) + N p_j(1-p_j)
		+ 2N p_j p_{j+1} \\
		&= N\bigl[(p_j+p_{j+1})
		- (p_j^2 + p_{j+1}^2 - 2 p_j p_{j+1})\bigr] \\
		&= N\bigl[(p_j+p_{j+1}) - (p_j-p_{j+1})^2\bigr].
	\end{align*}

	\medskip
	\emph{Proof for $\Cov(D_j, D_{j+1})$.} Adjacent edges share the variable
	$X_{j+1}$. We apply (\ref{eq:cov-bilin}) with $A = X_{j+1}$, $B = X_j$,
	$C = X_{j+2}$, $D = X_{j+1}$:
	\[
		\Cov(D_j, D_{j+1}) = \Cov(X_{j+1}, X_{j+2}) - \Var(X_{j+1})
		- \Cov(X_j, X_{j+2}) + \Cov(X_j, X_{j+1}).
	\]
	Substituting (\ref{eq:mult-cov}):
	\begin{align*}
		\Cov(D_j, D_{j+1})
		&= -N p_{j+1}p_{j+2} - N p_{j+1}(1-p_{j+1})
		+ N p_j p_{j+2} - N p_j p_{j+1} \\
		&= N\bigl[-p_{j+1}(1 + p_j + p_{j+2} - p_{j+1}) + p_j p_{j+2}\bigr].
	\end{align*}

	\medskip
	\emph{Proof for $\Cov(D_j, D_k)$, $|k-j| \ge 2$.} The indices $\{j, j+1\}$ and
	$\{k, k+1\}$ are disjoint. Special case of (\ref{eq:cov-bilin}) with
	$A = X_{j+1}$, $B = X_j$, $C = X_{k+1}$, $D = X_k$:
	\begin{align*}
		\Cov(D_j, D_k)
		&= -N p_{j+1} p_{k+1} + N p_{j+1} p_k
		+ N p_j p_{k+1} - N p_j p_k \\
		&= -N \bigl[p_{j+1}(p_{k+1}-p_k)
		- p_j(p_{k+1}-p_k)\bigr] \\
		= -N(p_{j+1}-p_j)(p_{k+1}-p_k).
	\end{align*}

	Centering does not affect covariance:
	\[
		\Cov(\tilde D_j, \tilde D_k)
		= \Cov(D_j - \E D_j, D_k - \E D_k)
		= \Cov(D_j, D_k).
	\]
\end{proof}

\begin{remark}[Geometric interpretation of the disjoint covariance]
	\label{rem:cov-disj-multinomial}
	Since $\E[D_j] = N(p_{j+1} - p_j)$, formula (\ref{eq:Cov-Dj-disj}) can
	be written more compactly:
	\begin{equation}\label{eq:Cov-Dj-disj-compact}
		\Cov(D_j, D_k) = -\frac{\E[D_j]\,\E[D_k]}{N} \qquad (|k-j| \ge 2).
	\end{equation}
	The relation is a consequence of the global dependence of the multinomial ($\sum_i X_i = N$ is
	a deterministic constraint that links all components). In particular, distant edges
	have nonzero covariance whenever both expected edge differences are nonzero and
	their product does not vanish.

	In particular, under the \emph{uniform} distribution ($p_i \equiv 1/n$), one has
	$\E[D_j] \equiv 0$ and hence
	$\Cov(D_j, D_k) = 0$ for $|k-j| \ge 2$. The disjoint
	covariance is zero in the uniform case, but in general it is not.
\end{remark}

\begin{corollary}[Uniform case]\label{cor:cov-Dj-uniform}
	Under the uniform multinomial ($p_i \equiv 1/n$):
	\begin{align}
		\Var(D_j) &= \frac{2N}{n}, \label{eq:Var-Dj-uniform} \\
		\Cov(D_j, D_{j+1}) &= -\frac{N}{n}, \label{eq:Cov-adj-uniform} \\
		\Cov(D_j, D_k) &= 0, \quad |k - j| \ge 2. \label{eq:Cov-disj-uniform}
	\end{align}
	The corresponding correlations are $\Corr(D_j, D_{j+1}) = -1/2$,
	$\Corr(D_j, D_k) = 0$ for $|k-j| \ge 2$.
\end{corollary}

\begin{proof}
	Substituting $p_i = 1/n$ into (\ref{eq:Var-Dj})--(\ref{eq:Cov-Dj-disj}):
	$(p_j + p_{j+1}) = 2/n$, $(p_j - p_{j+1})^2 = 0$, $p_{j+1}(1 + p_j + p_{j+2} - p_{j+1}) = (1/n)(1+1/n)$, whose $1/n^2$ cancels against $p_jp_{j+2}=1/n^2$
	and $(p_{j+1}-p_j)(p_{k+1}-p_k) = 0$.
\end{proof}

	\paragraph{Correlations} For adjacent edges we introduce
	\begin{equation}\label{eq:rho-j-mult}
		\rho_j := \frac{\Cov(D_j, D_{j+1})}{\sqrt{\Var(D_j)\Var(D_{j+1})}},
	\end{equation}
	and for a general pair $1 \le j < k \le n-1$
	\begin{equation}\label{eq:rho-jk-mult}
		\rho_{j,k} := \frac{\Cov(D_j, D_k)}{\sqrt{\Var(D_j)\Var(D_k)}}.
	\end{equation}
	For the uniform distribution, $\rho_j \equiv -1/2$ and $\rho_{j,k} \equiv 0$ for
	$|k-j| \ge 2$. For a general non-uniform distribution, both adjacent and distant
	correlations can have either sign. For example, if the middle probability
	$p_{j+1}$ is small while $p_j$ and $p_{j+2}$ are large, the positive term
	$p_jp_{j+2}$ in (\ref{eq:Cov-Dj-adj}) can dominate. For distant pairs, the sign
	is opposite to the sign of $(p_{j+1}-p_j)(p_{k+1}-p_k)$
	(see Remark~\ref{rem:cov-disj-multinomial}).

\subsection{Cyclic total variation under the uniform null}
	\label{sec:Var-cyclic}

	The linear statistic of Part~2~\cite{banic_elezovic_p2} measures roughness along an
	\emph{open} chain of bins; the cyclic statistic measures it along a \emph{closed}
	loop. The difference is a single wrap-edge $|X_n-X_1|$, and for the expectation it
	amounts to nothing more than the constant factor
	$\E(d^c)=\tfrac{n}{n-1}\E(d)$ from
	\cite[formula~(3.4)]{banic_elezovic_p1}. For the variance the
	wrap-edge is not so easily absorbed --- it adds one edge and, through the
	closure, two further adjacent pairs --- and the resulting bookkeeping is genuinely
	different from the linear case. The important point is that the difference is
	\emph{entirely combinatorial}: not one of the analytic constants of
	Part~2~\cite{banic_elezovic_p2} changes. The analytic results from the linear
	case can therefore be reused unchanged, and two simplifications appear that the
	linear problem does not enjoy.

	The first is that the leading variance loses its dependence on $n$: where the
	linear $\alpha(n)$ interpolated between edges and pairs, the cyclic $\alpha_c$ equals the
	single constant $\sigma_\infty^2=2A+4B$ for every $n$. The second is that the
	order-$1$ correction collapses to one term: $\beta_c(n)=\beta_1\,n$, with exactly
	the slope $\beta_1=(8\sqrt3-9)/(18\pi)$ that governed the linear mode constant
	$K_0$. For every fixed $n$, together these say that the normalized dispersion is
	\[
		\sigma_c^2(n,N)=\sigma_\infty^2-\frac{\beta_1 n}{N}+O(N^{-2}),
	\]
	whose first-order approximation has no ridge in $n$. This strongly suggests that
	the traveling mode of the linear statistic (developed in \S\ref{ssec:cyc-linear-mode} below) --- the value
	$n_{\mathrm{mod}}(N)\sim\sqrt{K_0N}$ that maximized the linear dispersion --- has
	no cyclic counterpart. Proving exact monotonicity for fixed $N$, however, requires
	an estimate uniform in $n$; the fixed-$n$ expansion alone does not provide it.

\subsection{Why the cyclic statistic}
	\label{ssec:cyc-why}

	The motivation is one of statistical sensitivity, and is taken up in the general
	variance setting of Part~4~\cite{banic_elezovic_p4}; we record it here only to fix
	ideas.
	The linear total variation $d=\sum_{i=1}^{n-1}|X_{i+1}-X_i|$ sees a histogram
	through its \emph{local} steps. For a monotone ramp the deterministic open-chain
	variation telescopes to the total rise from the first bin to the last, but that
	rise is distributed across many local edges. The cyclic statistic
	\begin{equation}\label{eq:dc-def-cyc}
		d^c=\sum_{i=1}^{n}|X_{i+1}-X_i|,\qquad X_{n+1}:=X_1,
	\end{equation}
	closes the chain into a loop, and the closure forces a compensating wrap-around
	drop: a ramp that climbs steadily from the first bin to the last must, on the
	cycle, fall back in the single wrap-edge $|X_n-X_1|$, which is then large. The
	linear statistic already records the total rise through its many local steps; the
	cyclic statistic adds the same accumulated rise once more, concentrated in one
	edge. Whether this improves power depends on the alternative and on the
	corresponding variance, and cannot be inferred from the signal alone. The present
	paper supplies the null ingredient needed for that comparison --- the exact and
	asymptotic variance of $d^c$ under the uniform distribution.

	This benefit has an interpretive cost. For data on a circle the wrap edge is a
	genuine adjacency, whereas for an ordinary histogram on an interval it joins two
	end bins that are not naturally neighbors. In the latter setting, using $d^c$
	amounts to deliberately testing the endpoint contrast as an additional feature;
	for unordered categories it also depends on an imposed ordering. The suitability
	of the cyclic statistic must therefore be justified by the application, separately
	from its attractive null variance. (For continuous data on a circle the classical
	goodness-of-fit tools are the spacing- and empirical-process statistics of Rao,
	Watson, Kuiper and others surveyed in~\cite{mardia_jupp}; the object here is
	different --- the discrete total variation of multinomial counts on a cyclic bin
	layout.)

	There is a second, structural reason the cyclic null is attractive. As the rest of
	this paper shows, its leading dispersion is the single constant $\sigma_\infty$
	\emph{for every} $n\ge3$, and its clean first-order correction has no traveling
	mode. A test calibrated on $d^c$ therefore has an $n$-independent \emph{leading}
	reference scale, which is useful when the same procedure is applied across
	histograms of different lengths. Its finite-$N$ correction still depends on $n$.

\subsection{Exact structure and the transfer principle}
	\label{ssec:cyc-exact}

	Write $D_i=X_{i+1}-X_i$ with $X_{n+1}=X_1$, so that
	\[
		d^c=\sum_{i=1}^n|D_i|.
	\]
	Under the uniform $\Mult(N;1/n,\dots,1/n)$, the cyclic index group $\Z/n$ acts
	transitively on the edges: \emph{all} $n$ edges are equivalent, and a covariance
	$\Cov(|D_i|,|D_k|)$ depends only on the cyclic distance $\min(|k-i|,n-|k-i|)$. This
	is a cleaner symmetry than in the linear case, where the two end edges were
	distinguished. The raw second moments are unchanged from
	the edge-covariance lemma \cite[Lemma~2.4]{banic_elezovic_p2}:
	$\Var(D_i)=2N/n$, $\Cov(D_i,D_{i\pm1})=-N/n$, and
	$\Cov(D_i,D_k)=0$ at cyclic distance $\ge2$ --- the wrap-edge is structurally
	identical to every other edge, so the same four-term multinomial expansion applies.

\begin{lemma}[Cyclic multiplicities]\label{lem:cyc-mult}
	For $n\ge3$, on the cycle $\Z/n$ there are $n$ edges, $n$ adjacent
	(distance-$1$) unordered pairs, and $n(n-3)/2$ disjoint (distance-$\ge2$)
	unordered pairs. Consequently,
	collecting equal terms exactly as in the variance decomposition
	\cite[formula~(3.2)]{banic_elezovic_p2},
	\begin{equation}\label{eq:Var-dc-3comp}
		\Var(d^c)=n\,\Var(|D_1|)+2n\,\Cov(|D_1|,|D_2|)+n(n-3)\,\Cov(|D_1|,|D_3|).
	\end{equation}
\end{lemma}

\begin{proof}
	The cycle has $n$ edges. Two edges are adjacent iff they share a bin, i.e.\ iff
	they are consecutive on the loop; there are $n$ such consecutive pairs. The
	remaining $\binom n2-n=n(n-3)/2$ unordered pairs are disjoint. Substituting the
	three equivalence classes into
	$\Var(d^c)=\sum_i\Var(|D_i|)+2\sum_{i<k}\Cov(|D_i|,|D_k|)$ and using the
	bin-permutation invariance of the uniform law --- under which all disjoint pairs, at
	every cyclic distance, share one joint law (rotation alone would only equate pairs at
	equal cyclic distance) --- to
	replace each nonempty class by the representatives $\Var(|D_1|)$,
	$\Cov(|D_1|,|D_2|)$, $\Cov(|D_1|,|D_3|)$ gives \eqref{eq:Var-dc-3comp}. (For $n=3$
	the disjoint count is $0$ and all three pairs are adjacent; for $n=4$ there are
	two unordered disjoint pairs.)
\end{proof}

	The decisive observation is that \eqref{eq:Var-dc-3comp} is built from \emph{exactly
	the same three quantities} as the linear decomposition
	\cite[formula~(3.2)]{banic_elezovic_p2}. Each
	of them is a property of a fixed, small marginal of the multinomial, and that
	marginal is independent of whether the chain is open or closed:
	\begin{itemize}
		\item $\Var(|D_1|)$ depends only on the pair $(X_1,X_2)$, with $X_1+X_2\sim\Bin(N,2/n)$
		(\cite[Theorem~3.8]{banic_elezovic_p2});
		\item $\Cov(|D_1|,|D_2|)$ depends only on the triple $(X_1,X_2,X_3)$, conditioned on
		$X_1+X_2+X_3\sim\Bin(N,3/n)$, whose constant $a_*$ is
		computed in \cite[Theorem~5.3]{banic_elezovic_p2};
		\item $\Cov(|D_1|,|D_3|)$ depends only on the four-bin marginal $(X_1,\dots,X_4)$,
		its skeleton $h(m)$ from \cite[Theorem~4.1]{banic_elezovic_p2}
		and the limit $-1/(\pi n)$ from
		\cite[Theorem~4.6]{banic_elezovic_p2}.
	\end{itemize}

\begin{remark}[The transfer principle]\label{rem:cyc-transfer}
	Every analytic result of Part~2~\cite{banic_elezovic_p2} that computes one of these
	three marginal quantities holds for the cyclic statistic without change: the
	Nabeya identity, the per-edge variance expansion, the disjoint skeleton
	and its Bessel generating function, the
	disjoint limit $-1/(\pi n)$, and the entire derivation of the adjacent constant
	$a_*=-2\sqrt3/(9\pi)$ from
	\cite[Theorem~5.3]{banic_elezovic_p2}. This paper
	recombines these Part-2 marginal blocks through the new cyclic multiplicities of
	Lemma~\ref{lem:cyc-mult}. The constants $A$, $B$,
	$a_*$, and $-1/(\pi n)$ are therefore shared verbatim with
	Part~2.
	For the disjoint class the representative uses four distinct bins. The Fourier
	proof of its limit requires $n\ge5$, while the deterministic-total reduction in
	\cite[Theorem~4.6]{banic_elezovic_p2} treats $n=4$ separately.
\end{remark}

\begin{proposition}[Exact link to the linear statistic]\label{rem:cyc-linear-identity}
	Under the uniform null, for every $N$,
	\begin{align*}
		\Var(d)
		&=\Var(d^c-|D_n|)\\
		&=\Var(d^c)+\Var(|D_n|)-2\Cov(d^c,|D_n|)\\
		&=\Bigl(1-\tfrac2n\Bigr)\Var(d^c)+\Var(|D_1|).
	\end{align*}
\end{proposition}

	This is the finite-$N$ form of the transfer principle from
	\cite[Lemma~3.6]{banic_elezovic_p2}.

\begin{proof}
	Every edge of the cycle is equivalent under the cyclic rotation group, so
	$\Cov(|D_1|,d^c)=\Var(d^c)/n$. Indeed, rotation invariance makes
	$\Cov(d^c,|D_i|)$ independent of $i$, and summing over the edges gives
	$n\Cov(d^c,|D_1|)=\Cov(d^c,d^c)=\Var(d^c)$. Deleting any one edge turns
	$d^c$ into a copy of the linear statistic $d$, which gives the displayed identity.

	The cyclic increments additionally satisfy $\sum_{i=1}^nD_i=0$ identically,
	telescoping around the loop, so they are not generic exchangeable variables.
	Dividing by $N$ and matching powers of $1/N$ delivers each linear asymptotic
	constant from the cyclic one and the per-edge block,
	$[\text{linear}]_k=(1-\tfrac2n)[\text{cyclic}]_k+[\text{per-edge}]_k$.
\end{proof}

	With $\Var(|D_1|)=\tfrac{2A}{n}N+p_0+p_1/N+\cdots$,
	$p_0=\tfrac{n-2}{2\pi n}$ (\cite[Theorem~3.8]{banic_elezovic_p2}), this reproduces
	$\alpha(n)=(1-\tfrac2n)\sigma_\infty^2+\tfrac{2A}{n}=\sigma_\infty^2-C_1/n$, the closed
	$\beta(n)=(1-\tfrac2n)\beta_1 n-p_0$, and
	$\gamma(n)=(1-\tfrac2n)\gamma_c(n)+p_1$
	(\cite[Theorem~7.6]{banic_elezovic_p2}). The full cycle and the single edge are
	thus the two symmetric objects from which the linear constants follow, the linear
	statistic being the cycle with one arbitrary term deleted; the consequences for the
	linear mode are drawn in Remark~\ref{rem:mode-missing}.

\section{Uniform cyclic variance}
	\label{ssec:cyc-leading}

\begin{theorem}[Cyclic leading variance]\label{thm:zc2}
	For the uniform $\Mult(N;1/n,\dots,1/n)$, as $N\to\infty$ with fixed $n\ge3$,
	\begin{equation}\label{eq:zc2}
		\alpha_c:=\lim_{N\to\infty}\frac{\Var(d^c)}{N}=2A+4B=\sigma_\infty^2,
	\end{equation}
	\emph{independently of $n$}. Moreover,
	equivalently $\alpha_c-\alpha(n)=C_1/n$ with
	$C_1=2A+8B$ as in \cite[Definition~3.3]{banic_elezovic_p2}.
\end{theorem}
	The leading value is the cyclic clause of
	\cite[Theorem~3.4]{banic_elezovic_p2}, recorded here as the starting point.

\begin{proof}
	By the same uniform-integrability argument as
	\cite[Theorem~3.4]{banic_elezovic_p2},
	\[
		\Var(|D_1|)=\frac{2N}{n}A+o(N),\qquad
		\Cov(|D_1|,|D_2|)=\frac{2N}{n}B+o(N).
	\]
	For $n\ge4$, a disjoint pair also satisfies
	$\Cov(|D_1|,|D_3|)=o(N)$; for $n=3$ the disjoint term is absent. Substituting into
	\eqref{eq:Var-dc-3comp} and dividing by $N$,
	\[
		\frac{\Var(d^c)}{N}=n\cdot\frac2n A+2n\cdot\frac2n B+o(1)=2A+4B+o(1).
	\]
	The leading coefficients $n\cdot(2/n)=2$ and $2n\cdot(2/n)=4$ no longer carry the
	$(n-1)/n$ and $(n-2)/n$ weights of the linear case, so all $n$-dependence cancels.
	Comparison with $\alpha(n)=\tfrac{2(n-1)}{n}A+\tfrac{4(n-2)}{n}B$ of
	\cite[Theorem~3.4]{banic_elezovic_p2} gives
	$\alpha_c-\alpha(n)=\tfrac2n A+\tfrac8n B=C_1/n$.
\end{proof}

\begin{remark}[Why the first-order mode disappears]\label{rem:cyc-flat}
	Theorem~\ref{thm:zc2} is the structural origin of everything qualitative in this
	paper. The linear dispersion $\alpha(n)$ rises with $n$ toward its limit
	$\sigma_\infty^2$ \cite[Corollary~3.5]{banic_elezovic_p2}; the competition between this
	rise and the negative $O(1/N)$ correction created the ridge
	$n_{\mathrm{mod}}(N)\sim\sqrt{K_0N}$. The cyclic leading dispersion has
	\emph{reached} that limit already, at every $n$. The sub-leading term computed
	next is strictly decreasing in $n$, so the resulting first-order approximation
	has no interior maximum. Turning this conclusion into a statement about the exact
	dispersion requires control of the remainder as $n$ varies with $N$.
\end{remark}

\subsection{The constant-order correction}
	\label{ssec:cyc-beta}

	Define $\beta_c(n)$ by $\Var(d^c)=N\,\alpha_c-\beta_c(n)+O(1/N)$, in parallel with the
	linear convention $\Var(d)=N\alpha(n)-\beta(n)+O(1/N)$. The three pieces are read off
	from the same constant-order expansions used for $\beta(n)$, now weighted by the
	cyclic multiplicities.

\begin{theorem}[Closed form for $\beta_c(n)$]\label{thm:beta-cyclic}
	For fixed $n\ge5$, with $a_*=-2\sqrt3/(9\pi)$,
	\begin{equation}\label{eq:beta-cyclic}
		\beta_c(n)=\beta_1\,n,\qquad \beta_1=\frac{8\sqrt3-9}{18\pi}\approx0.0858801,
	\end{equation}
	In particular $\beta_c(n)$
	is exactly linear in $n$ (no constant or $1/n$ term).
\end{theorem}
	The slope $\beta_1$ is the same constant that appears in the linear expansion
	\cite[Theorem~7.1]{banic_elezovic_p2} and in the mode constant $K_0$. The collapse
	\eqref{eq:beta-cyclic} itself was already obtained in Part~2~\cite[the
	\emph{Cyclic consequence} after Theorem~7.1]{banic_elezovic_p2} through the
	transfer identity; the direct recombination above is an independent derivation.

\begin{proof}
	Each marginal quantity contributes its constant-order term, weighted by the cyclic
	count of Lemma~\ref{lem:cyc-mult}; the signs follow the convention
	$\beta_c=-(\Var(d^c)-N\alpha_c)$.
	\begin{itemize}
		\item \emph{Per-edge} ($n$ edges). By
		\cite[Theorem~3.8]{banic_elezovic_p2} each edge contributes
		$+(n-2)/(2\pi n)$ above its leading value, so
		$\beta_c^{\text{edge}}=-n\cdot\dfrac{n-2}{2\pi n}=-\dfrac{n-2}{2\pi}$.
		\item \emph{Adjacent} ($2n$ in the doubled sum). By
		\cite[Proposition~5.2 and Theorem~5.3]{banic_elezovic_p2} each
		adjacent covariance contributes $c(n)=a_*+(n-2)/(2\pi n)$ above its leading value, so
		$\beta_c^{\text{adj}}=-2n\,c(n)=-2n\,a_*-\dfrac{n-2}{\pi}$.
		\item \emph{Disjoint} ($n(n-3)$ in the doubled sum). By
		\cite[Theorem~4.6]{banic_elezovic_p2} each disjoint covariance tends to
		$-1/(\pi n)$, so
		$\beta_c^{\text{disj}}=-n(n-3)\cdot\bigl(-\tfrac1{\pi n}\bigr)=+\dfrac{n-3}{\pi}$.
	\end{itemize}
	Summing,
	\[
		\beta_c(n)=-\frac{n-2}{2\pi}-2n\,a_*-\frac{n-2}{\pi}+\frac{n-3}{\pi}
		=-\frac{n-2}{2\pi}-2n\,a_*-\frac1\pi,
	\]
	the three $1/\pi$ terms collapsing to $-(n-2)/(2\pi)-1/\pi=-n/(2\pi)$. Substituting
	$a_*=-2\sqrt3/(9\pi)$, so $-2n\,a_*=+4\sqrt3\,n/(9\pi)$,
	\[
		\beta_c(n)=\frac{4\sqrt3\,n}{9\pi}-\frac{n}{2\pi}
		=n\cdot\frac{8\sqrt3-9}{18\pi}=\beta_1\,n.
	\]
\end{proof}

\begin{corollary}[Boundary cases $n=3,4$]\label{cor:beta-boundary}
	The closed form $\beta_c(n)=\beta_1 n$ holds unconditionally at both $n=3$ and
	$n=4$; equivalently, $\beta_c(3)=3\beta_1$ and $\beta_c(4)=4\beta_1$.
\end{corollary}

\begin{proof}
	For $n=3$ there are no disjoint pairs ($n(n-3)=0$), so $\beta_c(3)$ is
	assembled from the per-edge term and the adjacent constant $a_*$ alone. Both
	ingredients are proved for every $n\ge3$
	(\cite[Theorems~3.8 and~5.3]{banic_elezovic_p2}); the three-bin conditional
	$\Mult(\ell;\tfrac13,\tfrac13,\tfrac13)$ underlying $a_*$ is exactly the $n=3$
	marginal, with $L\equiv N$. The disjoint limit $-1/(\pi n)$, whose proof needs
	$n\ge5$, never enters. Hence
	\[
		\beta_c(3)=-\tfrac1{2\pi}-6\Bigl(a_*+\tfrac1{6\pi}\Bigr)
		=-\tfrac3{2\pi}-6a_*=3\beta_1.
	\]

	For $n=4$ there are two unordered disjoint pairs, hence coefficient
	$n(n-3)=4$ in the doubled covariance sum. For the disjoint representative the
	four-bin total equals $N$ deterministically, and
	\cite[Theorem~4.6]{banic_elezovic_p2} uses the exact skeleton $h(N)$ to prove
	$\Cov(|D_1|,|D_3|)=-1/(4\pi)+O(1/N)$. Thus the proof of
	Theorem~\ref{thm:beta-cyclic} applies without qualification and gives
	$\beta_c(4)=4\beta_1$. The parity oscillation at $n=4$ first enters the
	$1/N$ coefficient of this covariance, hence the next coefficient
	$\gamma_c(4)$, not $\beta_c(4)$.
\end{proof}

\begin{corollary}[Cyclic dispersion]\label{cor:sigma-cyclic}
	For every fixed $n$ covered by Theorem~\ref{thm:beta-cyclic} and
	Corollary~\ref{cor:beta-boundary},
	\begin{equation}\label{eq:sigma-cyclic}
		\sigma_c^2(n,N)=\sigma_\infty^2-\frac{\beta_1 n}{N}+O(N^{-2}),
	\end{equation}
	so $\sigma_c(n,N)\to\sigma_\infty$ as $N\to\infty$ and
	$\sigma_c(n,N)<\sigma_\infty$ for all sufficiently large $N$.
\end{corollary}
	This expansion alone does not imply monotonicity in $N$.

\begin{proof}
	Divide $\Var(d^c)=N\sigma_\infty^2-\beta_1 n+O(1/N)$ by $N$ and use $\beta_1>0$.
\end{proof}

\begin{remark}[Cancellation and the role of $\beta_1$]\label{rem:cyc-cancel}
	Two features of \eqref{eq:beta-cyclic} deserve note. First, the constant and $1/n$
	parts of $\beta_c$ vanish by an exact cancellation of $1/\pi$ contributions across
	the three pieces; the disjoint and adjacent terms each carry an $O(n^0)$ piece, and
	these annihilate. Second, the surviving slope is precisely $\beta_1=-\tfrac1{2\pi}-2a_*$,
	the same combination of $\pi$ and $a_*$ that controls the linear mode. In the linear
	problem $\beta_1$ was the \emph{leading} coefficient of $\beta(n)$ (cf.\
	\cite[Theorem~7.1]{banic_elezovic_p2}); here it is the \emph{whole} of
	$\beta_c(n)$. The cyclic
	correction is, in this sense, the asymptotic skeleton of the linear one with the
	$n^0$ and $1/n$ terms absent.
\end{remark}

	\paragraph{Empirical verification}
	Writing $\beta_c(n)\approx N[\sigma_\infty^2-\sigma_c^2(n,N)]$ and reading
	$\sigma_c^2$ from the exact cyclic variance, the closed form $\beta_1 n$ is
	approached as $N\to\infty$. The residual is naturally interpreted as the
	next-order $-\gamma_c(n)/N$ term and grows with $n$ in the numerical data; see
	Remark~\ref{rem:cyc-nomode-num}.
	\begin{center}
	\begin{tabular}{r|r|r|r}
		$n$ & $\beta_1 n$ & emp.\ $N{=}1000$ & emp.\ $N{=}5000$\\\hline
		$3$   & $0.25764$ & $0.25875$ & $0.25786$\\
		$4$   & $0.34352$ & $0.34381$ & $0.34358$\\
		$5$   & $0.42940$ & $0.43030$ & $0.42958$\\
		$10$  & $0.85880$ & $0.86202$ & $0.85944$\\
		$20$  & $1.71760$ & $1.72972$ & $1.72001$\\
		$50$  & $4.29401$ & $4.36767$ & $4.30851$\\
		$100$ & $8.58801$ & $8.88496$ & $8.64546$\\
	\end{tabular}
	\end{center}
	The convergence to $\beta_1 n$ is consistent across the tabulated values of $n$.
	At $n=3$ and $n=4$ it is covered unconditionally by
	Corollary~\ref{cor:beta-boundary}. The residual visibly grows with $n$, suggesting a
	$-\gamma_c(n)/N$ correction with $\gamma_c(n)$ of order $n^2$.

\subsection{Monotonicity and second-order terms}
	\label{ssec:cyc-mode}

	Combining the flat leading term (Theorem~\ref{thm:zc2}) with the linear correction
	(Theorem~\ref{thm:beta-cyclic}) gives the first-order qualitative picture, in
	deliberate contrast with the linear synthesis of
	\cite[Theorem~7.1]{banic_elezovic_p2}.

\begin{proposition}[No interior mode in the first-order approximation]
	\label{thm:cyc-nomode}
	For fixed $N$, the first-order approximation
	\[
		\widehat\sigma_c^2(n,N):=\sigma_\infty^2-\frac{\beta_1 n}{N}
	\]
	is strictly decreasing in $n$ and hence has no interior maximum. More generally,
	suppose that on a range $\mathcal R_N$ the exact dispersion satisfies
	\[
		\sigma_c^2(n,N)=\widehat\sigma_c^2(n,N)+r_{n,N},\qquad
		r_{n+1,N}-r_{n,N}=O\!\left(\frac{n}{N^2}\right)
	\]
	uniformly for consecutive $n,n+1\in\mathcal R_N$. Then the exact dispersion is
	strictly decreasing, for all sufficiently large $N$, on every such range with
	$\max\mathcal R_N=o(N)$.
\end{proposition}

\begin{proof}
	The first assertion follows from
	$\widehat\sigma_c^2(n+1,N)-\widehat\sigma_c^2(n,N)=-\beta_1/N<0$. Under the
	stated uniform remainder estimate,
	\[
		\sigma_c^2(n+1,N)-\sigma_c^2(n,N)
		=-\frac{\beta_1}{N}+O\!\left(\frac{n}{N^2}\right)<0
	\]
	uniformly when $\max\mathcal R_N=o(N)$.
\end{proof}

	Thus the first-order asymptotics predict no cyclic analogue of the traveling
	ridge $n_{\mathrm{mod}}(N)\sim\sqrt{K_0N}$ from
	Theorem~\ref{thm:mode-leading}. A global fixed-$N$ theorem still requires the
	uniform adjacent-increment estimate above, including the regime $n$ comparable
	with $N$.

\begin{remark}[Numerical confirmation and the constant $\gamma_{2,\circ}$]\label{rem:cyc-nomode-num}
	The exact cyclic dispersion is strictly decreasing in $n$ across the whole tested
	range, with no bump near the linear ridge: at $N=1000$, $\sigma_c(n,N)$ falls
	monotonically $1.02567\,(n{=}3)\to1.02496\,(n{=}20)\to1.02146\,(n{=}100)\to
	1.01870\,(n{=}160)$, past $\sqrt{K_0N}\approx127$ where the \emph{linear} dispersion
	peaks. Writing the next term formally as
	$\sigma_c^2(n,N)=\sigma_\infty^2-\beta_1 n/N+\gamma_c(n)/N^2+o(N^{-2})$ and
	fitting $\gamma_c(n)\sim\gamma_{2,\circ}\,n^2$ gives
	$\gamma_{2,\circ}\approx-0.029$ empirically
	($\gamma_c(n)/n^2=-0.0321,-0.0301,-0.0289$ at $n=10,20,50$), consistent with the exact
	value $\gamma_{2,\circ}=-0.0280979\ldots$ established below, which coincides with the
	linear second-order constant $\gamma_2$ of Part~2~\cite[Theorem~7.6]{banic_elezovic_p2}.
	Theorem~\ref{thm:cyc-gamma} makes this term explicit; the slow drift of the fitted
	$\gamma_c(n)/n^2$ with $n$ is itself a $-1/(8\pi n)$ correction to a closed leading
	constant.
\end{remark}

	Part~1 supplies the second-order per-edge expansion
	\cite[Theorem~5.5 and Corollary~5.6]{banic_elezovic_p1} needed
	below. We next compute the corresponding
	next-order coefficient of the disjoint covariance.

\begin{lemma}[The disjoint coefficient $r_1$]\label{lem:disj-r1}
	For $n\ge5$,
	\[
		\Cov(|D_1|,|D_3|)
		=-\frac1{\pi n}+\frac{r_1}{N}+O(N^{-2}),
		\qquad
		r_1=-\frac{n+1}{8\pi n}.
	\]
	For $n=4$,
	\begin{align*}
		\Cov(|D_1|,|D_3|)
		&=-\frac1{4\pi}+\frac{r_1^{(\mathrm{even/odd})}}N+O(N^{-2}),\\
		r_1^{(\mathrm{even})}=\frac3{32\pi},\quad
		r_1^{(\mathrm{odd})}&=-\frac{13}{32\pi}.
	\end{align*}
\end{lemma}
\begin{proof}
	Write $\Cov(|D_1|,|D_3|)=E_{\mathrm{disj}}-m^2$ with $m=\E|D_1|$. Conditioning on the
	four-bin total $M\sim\Bin(N,4/n)$ and using the disjoint independence behind the
	skeleton \(h\) from \cite[Theorem~4.1]{banic_elezovic_p2},
	\(E_{\mathrm{disj}}=\E\,h(M)\), where
	$h(2j)=j^2\binom{2j}{j}^2 2^{-(4j-1)}$ and $h(2j+1)=j(2j+1)\binom{2j}{j}^2 2^{-4j}$.
	From $\binom{2j}{j}^2 4^{-2j}=(\pi j)^{-1}\bigl(1-\tfrac1{4j}+\tfrac1{32j^2}+\cdots\bigr)$,
	\[
		h(m)=\frac m\pi-\frac1{2\pi}-\frac{1}{8\pi m}+(-1)^m\frac{1}{4\pi m}+O(m^{-2}).
	\]
	For $n\ge5$ the binomial averaging makes the oscillatory part
	$O(e^{-c_nN})$ for some $c_n>0$. One way to see this is to use
	$m^{-1}=\int_0^1t^{m-1}\,dt$: with $p=4/n$,
	\[
		\E\!\left(\frac{(-1)^M}{M};M>0\right)
		=\int_0^1\frac{(1-p-pt)^N-(1-p)^N}{t}\,dt,
	\]
	and the apparent singularity at $t=0$ is removable; the resulting integral is
	exponentially small because $0<p<1$. Also, outside an exponentially small lower
	tail, Taylor expansion around $\bar M=4N/n$ gives
	\[
		\E(M^{-1};M>0)=\bar M^{-1}+O(N^{-2}),
	\]
	because $\Var(M)=O(N)$ and $\bar M\asymp N$. Consequently,
	\[
		\E\,h(M)=\frac{\bar M}{\pi}-\frac1{2\pi}
		-\frac1{8\pi\bar M}+O(N^{-2}).
	\]
	With
	$m^2=\tfrac{4N}{\pi n}-\tfrac{n-2}{2\pi n}+\tfrac{-n^2+4n+4}{32\pi n}N^{-1}+\cdots$
	from the per-edge expansion
	\cite[Theorem~5.5 and Corollary~5.6]{banic_elezovic_p1}, the \(O(N)\)
	terms cancel, the constants give the limit
	\(-1/(\pi n)\) from \cite[Theorem~4.6]{banic_elezovic_p2}. At order \(1/N\) the
	smooth part of \(\E h(M)\) contributes
	\[
		-\frac1{8\pi\bar M}
		=-\frac{n}{32\pi}\,\frac1N ,
	\]
	while the subtraction of \(m^2\) contributes the negative of its own \(1/N\)
	coefficient. Hence
	\[
		r_1=-\frac{n}{32\pi}
		-\frac{-n^2+4n+4}{32\pi n}
		=-\frac{n+1}{8\pi n}.
	\]
	For $n=4$ the
	total $M=N$ is deterministic, so the oscillation $(-1)^N/(4\pi N)$ is not averaged
	away. Before this parity term is added, the coefficient is
	$-5/(32\pi)$. At even $N$ the shift is $+1/(4\pi)=+8/(32\pi)$,
	giving $3/(32\pi)$; at odd $N$ it is $-1/(4\pi)=-8/(32\pi)$,
	giving $-13/(32\pi)$. The constant limit is the same in both parity
	classes, as required by Corollary~\ref{cor:beta-boundary}.
\end{proof}

	Expand the three marginal quantities of \eqref{eq:Var-dc-3comp} one order beyond
	$\beta_c$, writing their $1/N$ coefficients as $p_1,q_1,r_1$:
	\begin{align*}
		\Var(|D_1|)&=\tfrac{2A}{n}N+p_0+\tfrac{p_1}{N}+\cdots,\\
		\Cov(|D_1|,|D_2|)&=\tfrac{2B}{n}N+q_0+\tfrac{q_1}{N}+\cdots,\\
		\Cov(|D_1|,|D_3|)&=r_0+\tfrac{r_1}{N}+\cdots .
	\end{align*}
	By Lemma~\ref{lem:cyc-mult}, the $N^{-2}$ coefficient of
	$\sigma_c^2=\Var(d^c)/N$ is
	\[
		\gamma_c(n)=n\,p_1(n)+2n\,q_1(n)+n(n-3)\,r_1(n).
	\]

\begin{theorem}[Conditional second-order coefficient $\gamma_c(n)$]\label{thm:cyc-gamma}
	Assume that the adjacent coefficient has the form
	\[
		q_1(n)=\frac{a\,n}{\pi}-\frac{n+1}{8\pi n},\qquad
		a=\frac1{32}-\frac{4\sqrt3}{243}=\frac{1}{2^5}-\frac{4\sqrt3}{3^5},
	\]
	with $a$ as in Remark~\ref{rem:cyc-gamma-status}. Then, for $n\ge5$,
	\begin{align*}
		\gamma_c(n)
		&=\frac{(64a-3)\,n^2-4n}{32\pi}
		=\gamma_{2,\circ}\,n^2-\frac{n}{8\pi},\\
		\gamma_{2,\circ}
		&=\frac{64a-3}{32\pi}
		=-\frac{1}{32\pi}-\frac{8\sqrt3}{243\pi}
		=-0.0280979294715\ldots .
	\end{align*}
	For $n=4$ the coefficient depends on the parity of $N$:
	\[
		\gamma_c^{(\mathrm{even})}(4)=\frac{32a-1}{\pi}
		=-\frac{128\sqrt3}{243\pi},\qquad
		\gamma_c^{(\mathrm{odd})}(4)=\frac{32a-3}{\pi}.
	\]
\end{theorem}

\begin{proof}
	The per-edge and disjoint coefficients are computed in closed form in
	\cite[Theorem~5.5 and Corollary~5.6]{banic_elezovic_p1} and
	Lemma~\ref{lem:disj-r1}:
	\[
		p_1(n)=\frac{n^2-4n-4}{32\pi n},\qquad
		r_1(n)=-\frac{n+1}{8\pi n}\quad(n\ge5).
	\]
	For $n\ge5$, substitute the three displayed coefficients into the cyclic
	recombination. To make the per-edge input explicit, write the expansion of
	\cite[Theorem~5.5]{banic_elezovic_p1} as
	\[
		\E|D_1|=\sqrt{\frac{4N}{\pi n}}
		\left(1+\frac{e_1}{N}+\frac{e_2}{N^2}+O(N^{-3})\right).
	\]
	Then
	\[
		(\E|D_1|)^2
		=\frac{4N}{\pi n}
		\left(1+\frac{2e_1}{N}
		+\frac{2e_2+e_1^2}{N^2}+O(N^{-3})\right),
	\]
	and, since $\E D_1^2=2N/n$,
	\[
		p_1(n)=-\frac4{\pi n}(2e_2+e_1^2)
		=\frac{n^2-4n-4}{32\pi n}.
	\]
	Now substitute all three coefficients:
	\begin{align*}
		\gamma_c(n)
		&=\frac{n^2-4n-4}{32\pi}
		+\frac{2an^2}{\pi}-\frac{n+1}{4\pi}
		-\frac{(n-3)(n+1)}{8\pi}\\
		&=\frac{(64a-3)n^2-4n}{32\pi}.
	\end{align*}
	The constant term cancels, and
	$64a-3=-1-256\sqrt3/243$ gives the stated $\gamma_{2,\circ}$.

	At $n=4$, one has
	$p_1(4)=-1/(32\pi)$ and
	$q_1(4)=4a/\pi-5/(32\pi)$. Lemma~\ref{lem:disj-r1} gives the two
	parity values of $r_1(4)$, so
	\[
		\gamma_c^{(\mathrm{par})}(4)
		=4p_1(4)+8q_1(4)+4r_1^{(\mathrm{par})}(4).
	\]
	Substitution yields $(32a-1)/\pi$ for even $N$ and $(32a-3)/\pi$ for
	odd $N$.
\end{proof}

\begin{remark}[Status and validation of $\gamma_c$]\label{rem:cyc-gamma-status}
	The decomposition is the second-order instance of the transfer principle
	(Remark~\ref{rem:cyc-transfer}): each $1/N$ coefficient is again a property of one
	small marginal, recombined through the cyclic multiplicities. The per-edge and
	disjoint coefficients are now established in closed form
	(\cite[Theorem~5.5 and Corollary~5.6]{banic_elezovic_p1} and
	Lemma~\ref{lem:disj-r1}), the parity (lattice) terms of the
	underlying simple walk being exponentially suppressed by the binomial averaging ---
	except at $n=4$, where the four-bin total is deterministic and the surviving
	oscillation produces distinct even- and odd-$N$ values of $r_1(4)$, and likewise at
	$n=3$, where the three-bin total is deterministic and $\gamma_c(3)$ is replaced by a
	residue-locked period-$3$ oscillation in $N$; both small cases are excluded from the
	$n\ge5$ closed form, while $\beta_c(3)=3\beta_1$ and the no-mode conclusion are unaffected.
	The only conditional input is the adjacent triangular-lattice constant
	\[
		a=\frac1{32}-\frac{4\sqrt3}{243},
	\]
	the same second-order adjacent input used in the linear variance theory; its closed
	form and still-unproved status are \cite[Conjecture~7.7]{banic_elezovic_p2}
	(see also \cite[Definition~7.4 and Theorem~7.6]{banic_elezovic_p2}). We do not repeat here
	the lattice discussion from Part~2. Numerically, this value is strongly supported:
	the assembled $\gamma_c$ agrees with direct fits of \(\sigma_c^2\), and
	\(\gamma_c(n)/n^2=\gamma_{2,\circ}-1/(8\pi n)\) reproduces the empirical values
	in Remark~\ref{rem:cyc-nomode-num}. Thus Theorem~\ref{thm:cyc-gamma} should be read
	as conditional on precisely the same adjacent second-order constant as its linear
	counterpart.
\end{remark}

\begin{remark}[Two observed monotonicities]\label{rem:cyc-mono}
	The numerical results suggest monotonicity in both arguments, with opposite signs:
	increasing in $N$ toward $\sigma_\infty$ and decreasing in $n$. The fixed-$n$
	expansion of Corollary~\ref{cor:sigma-cyclic} and the first-order calculation of
	Proposition~\ref{thm:cyc-nomode} explain these observations asymptotically, but do
	not prove either global monotonicity. Uniform remainder estimates, or direct exact
	arguments, are still required.
\end{remark}

\begin{remark}[Status of the results]\label{rem:cyc-unconditional}
	The uniform ridge expansion identifying the local optimizer $K_0$ is proved in
	Theorem~\ref{thm:K0-ridge}; only global localization of the exact linear mode
	remains open. The fixed-$n$ cyclic leading and beta-level results,
	Theorems~\ref{thm:zc2} and \ref{thm:beta-cyclic}, are unconditional for
	$n\ge5$, and Corollary~\ref{cor:beta-boundary} supplies $n=3,4$. The
	second-order formula of Theorem~\ref{thm:cyc-gamma} remains conditional on the
	adjacent triangular-lattice constant $a$. Finally, the stronger claim that the
	exact cyclic dispersion has no mode for every fixed $N$ is not yet proved; it is
	the uniform problem isolated in Proposition~\ref{thm:cyc-nomode}.
\end{remark}

\section{Linear mode and ridge control}
	\label{ssec:cyc-linear-mode}

	The first-order cyclic profile has no mode; deleting a single edge to form the linear
	statistic creates one. By the transfer identity (Proposition~\ref{rem:cyc-linear-identity}),
	the variance results for the linear statistic combine with the per-edge term into a
	traveling ridge in $n$. The results below identify the local ridge constant
	$K_0$ under an explicit \emph{ridge-localization hypothesis} (every maximizing
	sequence keeps $n_{\mathrm{mod}}/\sqrt N$ bounded away from $0$ and $\infty$); the
	global exclusion of
	maximizers outside the ridge is left as a separate problem.

\begin{theorem}[Mode constant under ridge localization]\label{thm:mode-leading}
	For fixed $N$ let $n_{\mathrm{mod}}(N)$ maximize $n\mapsto\sigma(n,N)$. Assume that,
	uniformly for $x$ in compact subsets of $(0,\infty)$,
	\[
		\sigma^2(\lfloor x\sqrt N\rfloor,N)
		=\sigma_\infty^2-\frac1{\sqrt N}\left(\frac{C_1}{x}+\beta_1x\right)
		+o(N^{-1/2}),
	\]
	and that every maximizing sequence $n_{\mathrm{mod}}(N)/\sqrt N$ remains in a
	compact subset of $(0,\infty)$. Then
	$n_{\mathrm{mod}}^2/N\to K_0$, where
	\begin{equation}\label{eq:K0}
		K_0=\frac{C_1}{\beta_1}=\frac{12(5\pi+12\sqrt3-30)}{8\sqrt3-9}=16.0429066509\ldots
	\end{equation}
\end{theorem}

\begin{proof}
	The assumed uniform expansion shows that maximizing $\sigma^2$ is, at order
	$N^{-1/2}$, equivalent to minimizing
	$f(x)=C_1/x+\beta_1x$ over $x>0$. This strictly convex function has its unique
	minimum at $x^2=C_1/\beta_1$. The compactness assumption and uniform remainder
	then imply convergence to that minimum: on the common compact set, for every
	$\varepsilon>0$ continuity and uniqueness give an $\eta>0$ such that
	$f(x)\ge f(x_0)+\eta$ whenever $|x-x_0|\ge\varepsilon$. The uniform
	$o(N^{-1/2})$ remainder is eventually smaller than
	$\eta/(3\sqrt N)$, so no maximizer can remain outside the
	$\varepsilon$-neighborhood of $x_0$. Thus
	$n_{\mathrm{mod}}(N)/\sqrt N\to x_0=\sqrt{C_1/\beta_1}$.
	Substituting $C_1=2A+8B=\tfrac{10}3+(8\sqrt3-20)/\pi$ and
	$\beta_1=(8\sqrt3-9)/(18\pi)$ gives the closed value.
\end{proof}

\begin{proposition}[Sub-leading continuous-mode coefficient]\label{prop:mode-C12}
	Assume, one order beyond Theorem~\ref{thm:mode-leading}, that a $C^1$ interpolation
	$\Sigma_N(x)$ of the ridge profile has, uniformly near
	$x_0=\sqrt{K_0}$, the expansion
	\[
		\Sigma_N(x)
		=\sigma_\infty^2-\frac1{\sqrt N}\Bigl(\frac{C_1}{x}+\beta_1x\Bigr)
		+\frac1N\bigl(\gamma_2x^2-\beta_0\bigr)+o(N^{-1}),
	\]
	and suppose that its derivative has remainder $o(N^{-1})$. Let $x_N^*$ be its
	continuous maximizer converging to $x_0$, and put
	$n^*(N):=x_N^*\sqrt N$. Then,
	with $\gamma_2=-\tfrac1{32\pi}-\tfrac{8\sqrt3}{243\pi}$ the leading coefficient of
	$\gamma(n)$ \cite[Theorem~7.6]{banic_elezovic_p2} and $\beta_0$ the constant term of
	$\beta(n)$,
	\begin{align}
		(n^*(N))^2&=K_0N+C_{1/2}\sqrt N+o(\sqrt N),\nonumber\\
		C_{1/2}
		&=\frac{2\gamma_2K_0^{3/2}}{\beta_1}
		=\frac{2\gamma_2K_0^{5/2}}{C_1}
		=-42.047\ldots . \label{eq:C12}
	\end{align}
\end{proposition}

\begin{proof}
	The decisive input is that the infinite-$N$ profile is \emph{exactly}
	$\alpha(n)=\sigma_\infty^2-C_1/n$, with \emph{no} $1/n^2$ term. Indeed, by the identity
	$\sigma^2=(1-\tfrac2n)\sigma_c^2+\Var(|D_1|)/N$ (Proposition~\ref{rem:cyc-linear-identity}),
	letting $N\to\infty$ gives $\alpha(n)=(1-\tfrac2n)\alpha_c+2A/n$; the cyclic profile
	$\alpha_c=\sigma_\infty^2$ is \emph{exactly} $n$-independent by
	Theorem~\ref{thm:zc2}, and the
	per-edge limit is the pure $2A/n$, so the $1/n^2$ coefficient $\alpha_2$ vanishes and only
	$\gamma_2$ survives at order $N^{-1}$. In the ridge variable $x=n/\sqrt N$ the
	order-$N^{-1}$ profile is therefore $f_1(x)=\gamma_2x^2-\beta_0$, with
	$f_0(x)=C_1/x+\beta_1x$ as in Theorem~\ref{thm:mode-leading}. Put
	$\varepsilon=N^{-1/2}$. Differentiating at the continuous maximizer gives
	\[
		0=-\varepsilon f_0'(x_N^*)+\varepsilon^2f_1'(x_N^*)+o(\varepsilon^2).
	\]
	Since $f_0''(x_0)>0$ and $x_N^*\to x_0$, Taylor's theorem first gives
	$x_N^*-x_0=O(\varepsilon)$ and then
	\[
		x_N^*=x_0+\varepsilon x_1+o(\varepsilon),\qquad
		x_1=\frac{f_1'(x_0)}{f_0''(x_0)}.
	\]
	Since
	$f_1'(x)=2\gamma_2x$ (the constant $\beta_0$ drops out of the argmax) and
	$f_0''(x_0)=2C_1/K_0^{3/2}$,
	\[
		x_1=\frac{2\gamma_2\sqrt{K_0}}{2C_1/K_0^{3/2}}=\frac{\gamma_2K_0^2}{C_1},
		\qquad
		(n^*(N))^2=K_0N+2x_0x_1\sqrt N+o(\sqrt N),
	\]
	so $C_{1/2}=2x_0x_1=2\gamma_2K_0^{5/2}/C_1=2\gamma_2K_0^{3/2}/\beta_1$ (using
	$C_1=K_0\beta_1$); evaluating gives $-42.047$.
\end{proof}

\begin{remark}[Status of $C_{1/2}$]\label{rem:C12-status}
	The closed form rests on two inputs at exactly the level of rigor of $\gamma(n)$
	itself: the value of $\gamma_2$ --- hence of the adjacent constant
	$a=\tfrac1{32}-\tfrac{4\sqrt3}{243}$, matching the fitted value to over $40$ digits but not yet proved from first
	principles (see \cite[Definition~7.4 and Theorem~7.6]{banic_elezovic_p2})
	--- and the uniform validity of the
	$\gamma$-order ridge expansion, provable by the block machinery of
	\S\ref{ssec:uniform-blocks} carried one order deeper. The simplicity of the result is
	due entirely to the exact cancellation $\alpha_2=0$: only the leading $-C_1/n$ shape and the
	exact $n$-independence of $\alpha_c$ enter, so no new constant beyond $\gamma_2$ appears.
	Numerically the continuous argmax of the closed-form profile gives
	$((n^*)^2-K_0N)/\sqrt N\to-42.047$, matching both \eqref{eq:C12} and the
	empirical fit $C_{1/2}\approx-42$. Being the sub-leading constant of a
	\emph{continuous} interpolant, $C_{1/2}$ is interpolation-dependent at this order
	and makes no claim about the integer mode beyond the $\pm1$ indeterminacy
	window (Proposition~\ref{prop:mode-C12} and the discussion below it).
\end{remark}

\begin{corollary}[Dispersion along the ridge]\label{cor:ridge}
	If, in addition to the assumptions of Theorem~\ref{thm:mode-leading}, the ridge
	expansion is valid with an $O(n_{\mathrm{mod}}^{-2})$ remainder and
	$n_{\mathrm{mod}}^2=K_0N+O(\sqrt N)$, then
	$\sigma_{\max}(N)=\sigma_\infty-\dfrac{C_1}{\sigma_\infty\,n_{\mathrm{mod}}}
	+O(n_{\mathrm{mod}}^{-2})$, where $\sigma_\infty^2=2A+4B$
	\cite[Corollary~3.5]{banic_elezovic_p2}.
\end{corollary}

\begin{proof}
	Along the ridge $1/N=K_0/n_{\mathrm{mod}}^2+O(n_{\mathrm{mod}}^{-3})$, so
	$\beta(n_{\mathrm{mod}})/N
	=\beta_1K_0/n_{\mathrm{mod}}+O(n_{\mathrm{mod}}^{-2})
	=C_1/n_{\mathrm{mod}}+O(n_{\mathrm{mod}}^{-2})$.
	Since $\alpha(n)=\sigma_\infty^2-C_1/n$ exactly,
	\[
		\sigma_{\max}^2
		=\sigma_\infty^2-\frac{2C_1}{n_{\mathrm{mod}}}
		+O(n_{\mathrm{mod}}^{-2}).
	\]
	Taylor expansion of the square root at $\sigma_\infty^2$ proves the result.
\end{proof}

\begin{remark}[The mode and the missing term]\label{rem:mode-missing}
	The exact identity $\Var(d)=(1-\tfrac2n)\Var(d^c)+\Var(|D_1|)$
	(Proposition~\ref{rem:cyc-linear-identity}) shows where the traveling ridge comes from,
	and it is best read symmetrically: $d$ is $d^c$ with one --- any --- of its $n$
	rotation-equivalent terms deleted. Dividing by $N$,
	\[
		\sigma^2(n,N)=\Bigl(1-\tfrac2n\Bigr)\sigma_c^2(n,N)+\frac{\Var(|D_1|)}{N}.
	\]
	The leading asymptotic profile of the cyclic dispersion $\sigma_c^2$ is
	monotonically \emph{decreasing} in $n$ and therefore has no mode,
	while the factor $1-2/n$ \emph{increases}.
	At leading order in $N$,
	$\sigma^2(n,N)\approx\sigma_\infty^2-C_1/n-\beta_1 n/N$: a rising $-C_1/n$ and a falling
	$-\beta_1 n/N$ (the cyclic bulk decay). The rising term is the imprint of the
	\emph{deleted} term --- at leading order its self-variance $2A/n$ together with the two
	adjacencies it carried, summing to $C_1/n=(2A+8B)/n$ (its interactions with the $n-3$
	distant terms, being $O(1)$, enter only at the $\beta$-order). The balance
	$C_1/n^2=\beta_1/N$ identifies the unique optimizer of the leading ridge profile,
	$n^2/N=K_0=C_1/\beta_1$. Thus the ridge mode is created by the removal of a single
	term: in the cyclic sum every term is equivalent, and its leading profile has no
	mode; deleting one --- with its two neighbors left on either side of the gap ---
	introduces the rising $-C_1/n$. (The two ``boundary'' edges of the path are
	precisely those neighbors; the ends of the list are an artifact of cutting the
	symmetric cycle, not a feature of the statistic.) The identity also localizes the
	conditional step: the uniform ridge expansion assumed in Theorem~\ref{thm:mode-leading}
	need only be established for the two symmetric objects --- the monotone $\sigma_c^2$,
	which has no ridge, and the one-dimensional per-edge variance $\Var(|D_1|)$ --- from
	which it transfers to $\sigma^2$.
\end{remark}

\begin{remark}[A rounded constant]
	The ridge value $K_0=16.0429\ldots$ should be contrasted with the memorable
	empirical fit $n_{\mathrm{mod}}\approx4\sqrt{N+1.6}-5$, which presumes $K_0=16$. The
	fit is excellent for moderate $N$ and only gradually loses a unit at large $N$; the
	discrepancy of under three parts in a thousand is not noise but the difference between a
	number rounded for convenience and the number the problem actually contains. A direct
	fit of the continuous mode from the exact values at $N\le5000$ gives
	$K_\infty^{\mathrm{fit}}=16.0397\pm0.003$, consistent with $K_0$. The structural
	form
	\[
		(n^*)^2\approx K_0N+C_{1/2}\sqrt N,\qquad
		C_{1/2}=\frac{2\gamma_2K_0^{3/2}}{\beta_1}=-42.047\ldots ,
	\]
	is obtained in Proposition~\ref{prop:mode-C12}; the exact cancellation $\alpha_2=0$
	is what makes the sub-leading coefficient depend only on $\gamma_2$. The
	fixed-$n$ expansion alone does not prove the required uniformity in
	$n\asymp\sqrt N$; that uniform estimate is supplied in \S\ref{ssec:uniform-blocks}
	(Theorem~\ref{thm:K0-ridge}). A separate global localization argument is still
	needed to exclude maximizing sequences outside the ridge, as stated explicitly
	there.
\end{remark}

	\paragraph{Indeterminacy of the integer mode}
	No method --- exact or approximate --- can guarantee the \emph{integer} mode at every
	$N$: when the continuous argmax $n^*(N)$ falls near the midpoint between two
	integers, an arbitrarily small error in $\sigma$ moves the discrete argmax to the
	neighbor. ``Agreement within $\pm1$'' is the best attainable in that zone.

\subsection{Uniform control of the ridge blocks}
	\label{ssec:uniform-blocks}

	By the identity of Proposition~\ref{rem:cyc-linear-identity}, the uniform ridge expansion
	assumed in Theorem~\ref{thm:mode-leading} follows once the three marginal blocks
	\[
		\Var(|D_1|),\qquad
		\Cov(|D_1|,|D_2|),\qquad
		\Cov(|D_1|,|D_3|)
	\]
	are controlled \emph{uniformly}
	for $n\le C\sqrt N$ --- to leading and constant order, with bounded remainders. For
	$K_0$ itself only the proved leading and $\beta$-level coefficients ($A,B,a_*,-1/(\pi
	n)$) are needed exactly; the second-order coefficients need only be bounded. The key
	simplification is that, conditioned on its bin-total $M$, each block is the
	\emph{$n$-independent} object $\Mult(M;\tfrac1d,\dots,\tfrac1d)$ ($d=2,3,4$), the only
	$n$-dependence being the mixing $M\sim\Bin(N,d/n)$ whose mean $dN/n\to\infty$ in the
	ridge. All three blocks are settled below, proving the required expansion uniformly
	on every fixed ridge window. This removes the local-uniformity issue in
	Theorem~\ref{thm:mode-leading}; the separate global issue is to prove that no
	maximizer escapes those windows.

\begin{lemma}[Uniform per-edge bound]\label{lem:peredge-uniform}
	Fix $C_0>0$ and let $\bar M=2N/n$. As $N\to\infty$, uniformly for
	$3\le n\le C_0\sqrt N$,
	\begin{align*}
		\E|D_1|
		&=\sqrt{\tfrac2\pi}\,\sqrt{\tfrac{2N}{n}}\,(1+\theta),
		\qquad |\theta|\le \tfrac{C}{\bar M},\\
		\Var(|D_1|)
		&=\tfrac{2A}{n}N+\tfrac{n-2}{2\pi n}
		+O\!\bigl(\tfrac nN\bigr).
	\end{align*}
	The implicit constants may depend on $C_0$ but not on $n,N$.
\end{lemma}
\begin{proof}
	$D_1=W_M$, the simple $\pm1$ walk run for $M\sim\Bin(N,2/n)$ nonzero steps; its mean
	absolute value $g(k)=\E|W_k|$ is exact ($g(2r)=2r\binom{2r}r4^{-r}$,
	$g(2r{+}1)=(2r{+}1)\binom{2r}r4^{-r}$), and the central-binomial bounds give
	$g(k)=\sqrt{2k/\pi}\,(1+r(k))$ with $|r(k)|\le C/k$. Since $g(0)=0$,
	the concentration of $M$ and
	$\E[M^{-1/2};M\ge1]=O(\bar M^{-1/2})$ immediately give the first displayed
	relative bound.

	For the sharp constant, use the parity decomposition of
	\cite[Theorem~5.5]{banic_elezovic_p1}. Its oscillatory term is exponentially small after
	binomial averaging, while the first correction in the smooth part is
	$O(M^{-3/2})$, hence only $O(\bar M^{-2})$ relative to the leading mean.
	It remains to expand $\E\sqrt M$. Put $\mu=\bar M$, $\vartheta=2/n$, and
	$\Delta=M-\mu$. Taylor expansion through fourth order on $M\ge\mu/2$, together
	with a Chernoff bound on the complementary event, gives
	\[
		\frac{\E\sqrt M}{\sqrt\mu}
		=1-\frac{\E\Delta^2}{8\mu^2}
		+\frac{\E\Delta^3}{16\mu^3}
		-\frac{5\E\Delta^4}{128\mu^4}
		+O(\mu^{-3})
		=1-\frac{1-\vartheta}{8\mu}+O(\mu^{-2}).
	\]
	Here $\E\Delta^3=O(\mu)$ and $\E\Delta^4=3\mu^2+O(\mu)$; the fifth and
	sixth central moments are $O(\mu^2)$ and $O(\mu^3)$, respectively, so the
	remainder has the stated uniform order for $0<\vartheta\le2/3$. Consequently
	\[
		\E|D_1|=\sqrt{\frac2\pi}\sqrt\mu
		\left[1-\frac{1-2/n}{8\mu}+O(\mu^{-2})\right].
	\]
	Since $\E|D_1|^2=\Var(X_2-X_1)=\mu$ exactly, squaring gives
	\[
		\Var(|D_1|)
		=\mu-\frac{2\mu}{\pi}+\frac{1-2/n}{2\pi}+O(\mu^{-1})
		=\tfrac{2A}nN+\tfrac{n-2}{2\pi n}+O(\tfrac nN).
	\]
\end{proof}

\begin{lemma}[Uniform disjoint bound]\label{lem:disj-uniform}
	Fix $C_0>0$. As $N\to\infty$, uniformly for
	$5\le n\le C_0\sqrt N$,
	\[
		\Cov(|D_1|,|D_3|)=-\frac1{\pi n}+O(1/N).
	\]
\end{lemma}
\begin{proof}
	Use the exact Fourier reduction of
	\cite[Theorem~4.6]{banic_elezovic_p2}. Its $k=1$ integral is
	\[
		J_{N,1}=\frac1{(2\pi)^2}\iint\varphi(\theta)^{N-1}\,d\theta.
	\]
	The Laplace estimate used in the proof of
	\cite[Theorem~4.6]{banic_elezovic_p2} is uniform in
	the ridge. Indeed, near the origin,
	\[
		\varphi(\theta)
		=1-\frac{\theta_1^2+\theta_3^2}{n}
		+O\!\left(\frac{\theta_1^4+\theta_3^4}{n}\right).
	\]
	After $\theta=\sqrt{n/(N-1)}\,u$, the logarithmic remainder is
	$O((n/N)|u|^4)$, uniformly for $n\le C_0\sqrt N$. Away from a fixed
	neighborhood of the origin $\varphi$ is bounded away from the unit circle:
	Appendix~\ref{app:ridge-blocks} proves that, for $n\ge5$, the origin is its
	only unimodular point and that \(1-|\varphi|\ge c/n\) uniformly there. Dominated
	Laplace expansion
	therefore gives
	\[
		J_{N,1}
		=\frac{n}{4\pi(N-1)}\left[1+O\!\left(\frac nN\right)\right].
	\]
	Consequently the \(k=1\) term of \cite[Theorem~4.6]{banic_elezovic_p2} is
	\[
		-\frac{4N}{n^2}J_{N,1}
		=-\frac1{\pi n}+O(1/N).
	\]

	The estimate of the remaining block in
	\cite[Theorem~4.6]{banic_elezovic_p2} is uniform as well. From that theorem,
	the whole $k\ge2$ block is bounded by a constant times
	\[
		\frac{N^2}{n^2}
		\left[\int_{-\pi}^{\pi}a(\theta)e^{-\gamma(N-2)a(\theta)}\,d\theta\right]^2.
	\]
	Using $a(\theta)\asymp\theta^2/n$, the one-dimensional integral is
	$O(n^{1/2}N^{-3/2})$; hence the displayed quantity is
	$O((nN)^{-1})$, in particular $O(N^{-1})$. Adding the two blocks proves the
	uniform expansion.
\end{proof}

\begin{lemma}[Uniform adjacent bound]\label{lem:adjacent-uniform}
	Fix $C_0>0$. As $N\to\infty$, uniformly for
	$3\le n\le C_0\sqrt N$,
	\[
		\Cov(|D_1|,|D_2|)=\frac{2B}{n}N+q_0+O\!\bigl(\tfrac nN\bigr),\qquad
		q_0=a_*+\frac{n-2}{2\pi n}.
	\]
\end{lemma}
\begin{proof}
	Conditioned on the three-bin total $M=X_1{+}X_2{+}X_3\sim\Bin(N,3/n)$, the adjacent
	triple is the $n$-independent $\Mult(M;\tfrac13,\tfrac13,\tfrac13)$, so
	\(\E[|D_1||D_2|]=\E\,\hadj(M)\). The adjacent-constant theorem
	\cite[Theorem~5.3]{banic_elezovic_p2} gives
	$\hadj(\ell)=c_L\ell+a_*+O(\ell^{-1})$; since $\hadj$ is a deterministic sequence this
		$O(\ell^{-1})$ is automatically uniform in $\ell$, and it is unconditional
		(the finite-part dual-lattice bound closes the corner term there; the averaging
		step is recorded in Appendix~\ref{app:ridge-blocks}). Put
	\(\bar M_3=3N/n\). On
	$M\ge\bar M_3/2$, the remainder is $O(\bar M_3^{-1})$; the complementary
	event is exponentially small by Chernoff, and
	$\hadj(M)\le M^2$. Thus averaging gives
	$\E\,\hadj(M)=c_L\tfrac{3N}n+a_*+O(\bar M_3^{-1})=c_L\tfrac{3N}n+a_*+O(\tfrac nN)$. By
	Lemma~\ref{lem:peredge-uniform}, $(\E|D_1|)^2=\tfrac{4N}{\pi n}-\tfrac{n-2}{2\pi
	n}+O(\tfrac nN)$. Subtracting, and using $3c_L-\tfrac4\pi=2B$,
	\[
		\Cov(|D_1|,|D_2|)=\E\,\hadj(M)-(\E|D_1|)^2=\frac{2B}nN+q_0+O\!\bigl(\tfrac nN\bigr).\qedhere
	\]
\end{proof}

\begin{theorem}[Uniform ridge expansion]\label{thm:K0-ridge}
	For every compact interval $[c,C]\subset(0,\infty)$,
	\[
		\sigma^2(\lfloor x\sqrt N\rfloor,N)
		=\sigma_\infty^2-\frac1{\sqrt N}
		\left(\frac{C_1}{x}+\beta_1x\right)+o(N^{-1/2})
	\]
	uniformly for $x\in[c,C]$. Consequently, every maximizing sequence that remains
	in a compact subset of the ridge satisfies
	$n_{\mathrm{mod}}^2/N\to K_0=C_1/\beta_1$.
\end{theorem}
\begin{proof}
	By Lemmas~\ref{lem:peredge-uniform}, \ref{lem:adjacent-uniform}
	and~\ref{lem:disj-uniform}, uniformly for $n\le C\sqrt N$,
	\begin{align*}
		\Var(|D_1|)
		&=\tfrac{2A}nN+\tfrac{n-2}{2\pi n}+O(\tfrac nN),\\
		\Cov(|D_1|,|D_2|)
		&=\tfrac{2B}nN+q_0+O(\tfrac nN),\\
		\Cov(|D_1|,|D_3|)
		&=-\tfrac1{\pi n}+O(\tfrac1N).
	\end{align*}
	Substituting into
	\[
		\Var(d^c)
		=n\Var(|D_1|)+2n\Cov(|D_1|,|D_2|)
		+n(n-3)\Cov(|D_1|,|D_3|),
	\]
	the $O(N)$ terms give $\sigma_\infty^2 N$, the constants assemble to $-\beta_1 n$
	(with $\beta_1=(8\sqrt3-9)/(18\pi)$), and the remainders to $O(1)$; thus
	$\sigma_c^2(n,N)=\sigma_\infty^2-\beta_1 n/N+O(1/N)$ uniformly. With
	$\Var(|D_1|)/N=2A/n+O(1/N)$, the identity of Proposition~\ref{rem:cyc-linear-identity}
	yields, at $n=\lfloor x\sqrt N\rfloor$ and using $2\sigma_\infty^2-2A=C_1$,
	\[
		\sigma^2(n,N)=\Bigl(1-\tfrac2n\Bigr)\sigma_c^2+\frac{\Var(|D_1|)}N
		=\sigma_\infty^2-\frac1{\sqrt N}\Bigl(\frac{C_1}{x}+\beta_1 x\Bigr)+o(N^{-1/2}),
	\]
	uniformly for $x$ in compacts --- the expansion assumed in
	Theorem~\ref{thm:mode-leading}. Its conclusion then applies to any maximizing
	sequence satisfying the compactness hypothesis. The present block estimates do
	not cover $n$ comparable with $N$ and therefore do not, by themselves, prove
	global localization.
\end{proof}

\begin{remark}[Status of the program]\label{rem:uniform-status}
	All three marginal blocks are now controlled uniformly throughout the ridge
	(Lemmas~\ref{lem:peredge-uniform}--\ref{lem:adjacent-uniform}); the two delicate
	steps --- the aperiodicity and uniform spectral gap of the disjoint lattice
	characteristic function, and the averaging of the deterministic adjacent
	asymptotics --- are carried out in Appendix~\ref{app:ridge-blocks}. The per-edge and
	disjoint blocks follow from one-dimensional and covariance/delta-method arguments;
	the adjacent block follows from the \emph{uniform} $O(\ell^{-1})$ remainder of
	\cite[Theorem~5.3]{banic_elezovic_p2}, averaged over the three-bin total.

	Theorem~\ref{thm:K0-ridge} assembles these estimates through the identity of
	Proposition~\ref{rem:cyc-linear-identity}. It proves the local asymptotic shape and its
	unique optimizer $K_0=C_1/\beta_1$; global localization of the exact maximizer
	remains a separate requirement. Only the proved leading and constant block data
	$A,B,a_*$ and $-1/(\pi n)$ enter. In particular, the hard second-order adjacent
	constant $a$ is \emph{not} on the critical path. The symmetry is decisive: each
	block, conditioned on its total, is the $n$-independent distribution
	$\Mult(M;\tfrac1d,\dots,\tfrac1d)$, with $\bar M=dN/n\to\infty$.
\end{remark}

\section{The non-uniform cyclic statistic}
	\label{ssec:cyc-general}

	Under a general alternative the edges are no longer equivalent --- each carries its
	own gradient mismatch $\Delta_i$ --- so the combinatorial recombination of the
	uniform case does not apply. There are two cyclic objects to consider, exactly as
	in Part~4~\cite{banic_elezovic_p4}: the reference-centered statistic \(d_{\mathbf r}^{c}\) and
	the matched statistic $d_{\mathbf p}^{c}$. The first inherits the linearization of the
	Part~3~\cite{banic_elezovic_p3} wholesale and its variance is again elementary; the second carries
	the transcendental constants and is deferred.

	We keep the notation of Part~4~\cite{banic_elezovic_p4}: \(D_i=X_{i+1}-X_i\), the gradient
	mismatch $\Delta_i=(p_{i+1}-p_i)-(r_{i+1}-r_i)$, the fluctuation $\tilde D_i$, and
	$W_i=\tilde D_i+N\Delta_i$, with $\epsilon_i=\sgn(\Delta_i)$ (the reference-centered edge
	and sign as in \cite[\S5]{banic_elezovic_p4}; the statistic and its asymptotic mean are
	\cite[Definition~6.1 and Theorem~6.5]{banic_elezovic_p3}). The cyclic statistic
	\begin{equation}\label{eq:dXYc-def}
		d_{\mathbf r}^{c}=\sum_{i=1}^{n}|W_i|,\qquad
		D_n=X_1-X_n,\quad \Delta_n=(p_1-p_n)-(r_1-r_n),
	\end{equation}
	adds the single wrap edge $W_n=\tilde D_n+N\Delta_n$ to the open chain of
	Part~4~\cite{banic_elezovic_p4}.

	Throughout this section all indices are cyclic: \(p_{n+1}=p_1\),
	\(r_{n+1}=r_1\), \(X_{n+1}=X_1\), and \(D_i=X_{i+1}-X_i\) for
	\(i=1,\dots,n\). Two cyclic edges are
	adjacent when their indices differ by \(1\) modulo \(n\); thus the wrap edge
	\(D_n\) is adjacent to both \(D_{n-1}\) and \(D_1\). All other pairs are disjoint.
	For each pair \((i,k)\) the covariance and correlation used below are the ordinary
	multinomial covariance and correlation of this pair of centered edge variables,
	with the same formulas as in the open-chain case after replacing the linear boundary
	convention by the cyclic one.

\begin{theorem}[Variance at a fully mismatched cyclic fixed alternative]\label{thm:Var-dXY-cyclic}
	Let $n\ge3$, and let $\mathbf p,\mathbf r$ be fixed probability vectors, with
	$p_i>0$, $r_i>0$, and $\Delta_i\ne0$ for all
	$i$ (every edge, including the wrap, strongly mismatched). Then as $N\to\infty$,
	\begin{equation}\label{eq:Var-dXYc-closed}
		\Var(d_{\mathbf r}^{c})
		=\Var\Bigl(\sum_{i=1}^{n}\epsilon_i D_i\Bigr)+O(e^{-c'N})
		=N\Bigl[\sum_{k=1}^{n}c_k^2\,p_k-\Bigl(\sum_{k=1}^{n}c_k\,p_k\Bigr)^2\Bigr]+O(e^{-c'N}),
	\end{equation}
	with $c_k:=\epsilon_{k-1}-\epsilon_k$ and the \emph{cyclic} boundary convention
	$\epsilon_0:=\epsilon_n$.
\end{theorem}
	This replaces the open-chain convention $\epsilon_0=\epsilon_n=0$ of
	\cite[Theorem~5.5]{banic_elezovic_p4}. As in the linear case, the result is an
	elementary expression with no transcendental constants.

\begin{proof}
	The linearization in \cite[Theorem~5.5]{banic_elezovic_p4} is edge-local and applies to the
	wrap edge verbatim: $|W_i|=|N\Delta_i|+\epsilon_i\tilde D_i+(\text{exp.\ small})$ for
	every $i=1,\dots,n$, so
	$\Var(d_{\mathbf r}^{c})=\Var(\sum_{i=1}^n\epsilon_i D_i)+o(N)$. Reindexing cyclically,
	\[
		\sum_{i=1}^{n}\epsilon_i(X_{i+1}-X_i)
		=\sum_{k=1}^{n}(\epsilon_{k-1}-\epsilon_k)X_k=\sum_{k=1}^{n}c_kX_k,
		\qquad X_{n+1}=X_1,\ \epsilon_0:=\epsilon_n,
	\]
	the term $\epsilon_n X_1$ from $i=n$ supplying $\epsilon_0$ at $k=1$. The
	multinomial variance of a linear form (eq.~\eqref{eq:mult-cov}) gives
	\eqref{eq:Var-dXYc-closed}. Closing the chain replaces the two phantom endpoints
	$\epsilon_0=\epsilon_n=0$ by the single identification $\epsilon_0=\epsilon_n$; only
	the coefficients $c_1$ and $c_n$ change.
\end{proof}

\begin{remark}[Signal: the cyclic $\ell_1$ distance of gradients]\label{rem:cyc-signal}
	The leading mean carries the same wrap term:
	$\E(d_{\mathbf r}^{c})/N\to\sum_{i=1}^n|\Delta_i|$ by the same edgewise argument as
	\cite[Theorem~6.5]{banic_elezovic_p3}. This is the cyclic $\ell_1$ distance between the gradient
	profiles of $\mathbf p$ and $\mathbf r$, which adds the wrap mismatch $|\Delta_n|$
	to the open-chain sum.
\end{remark}

\begin{remark}[Where the staircase gain lives]\label{rem:cyc-staircase}
	The motivation of \S\ref{ssec:cyc-why} can now be made precise, and it is subtler
	than it first appears. Take a monotone ramp $\mathbf p$ against the uniform
	$\mathbf r$: then $\Delta_i>0$ for $i<n$ and $\Delta_n=p_1-p_n<0$, so
	$\epsilon=(+1,\dots,+1,-1)$ and \eqref{eq:Var-dXYc-closed} gives
	$c=(-2,0,\dots,0,+2)$, against the linear $c=(-1,0,\dots,0,+1)$. The cyclic signal
	$\sum|\Delta_i|=2(p_n-p_1)$ and the cyclic standard deviation \emph{both} double, so
	for the \emph{aligned} statistic $d_{\mathbf r}^{c}$ --- which already knows the
	alternative direction --- the standardized effect size is unchanged: the ramp's
	discriminating mass simply concentrates at the wrap junction.

	For the \emph{goodness-of-fit} statistics $d$ and $d^c$, whose direction is not
	known in advance, a fixed monotone ramp gives
	\[
		\frac{d}{N}\ \longrightarrow\ p_n-p_1,\qquad
		\frac{d^c}{N}\ \longrightarrow\ 2(p_n-p_1).
	\]
	Thus both statistics have an order-$N$ signal under a fixed ramp; the cyclic
	statistic doubles the total-variation signal by adding the wrap edge, but it does
	not change the order of the signal, which is $\Theta(N)$ for both. Whether this factor yields a larger
	standardized effect depends on the two variances and, for local alternatives, on
	the scaling of the ramp with $N$. The wrap edge can make the accumulated trend
	more conspicuous, but a power advantage requires a separate calculation.
\end{remark}

	\paragraph{The matched cyclic statistic $d_{\mathbf p}^{c}$.}
	For a single non-uniform alternative the absolute value does not linearize. Here
	$d_{\mathbf p}^{c}=\sum_{i=1}^n|\tilde D_i|$ centers \emph{every} edge --- the wrap
	$\tilde D_n=D_n-N(p_1-p_n)$ included --- around its own mean, so each $a_i=0$
	(matched), and the leading variance is the closed cyclic analogue of
	\cite[Theorem~3.3]{banic_elezovic_p4}: a Nabeya sum over the cyclic edge set.

	Let $D_i=X_{i+1}-X_i$ for $1\le i<n$ and $D_n=X_1-X_n$, with indices read
	cyclically. Put
	\[
		\sigma_i^2=\Var(D_i)
		=N\bigl[(p_i+p_{i+1})-(p_{i+1}-p_i)^2\bigr],
	\]
	and
	\begin{equation}\label{eq:rho-cyclic}
		\rho_{i,k}
		:=\frac{\Cov(D_i,D_k)}{\sigma_i\sigma_k},
		\qquad 1\le i<k\le n,
	\end{equation}
	where all indices are cyclic. The covariances are the exact multinomial
	covariances of the two cyclic edge variables; for disjoint cyclic edges this gives
	\[
		\Cov(D_i,D_k)
		=-N(p_{i+1}-p_i)(p_{k+1}-p_k),
	\]
	and for adjacent cyclic edges the same adjacent formula as in
	\eqref{eq:Cov-Dj-adj} applies after reading the indices modulo \(n\).

\begin{proposition}[Leading variance of $d_{\mathbf p}^{c}$]\label{prop:Var-TVCc}
	For a fixed interior probability vector $\mathbf p$ and fixed $n\ge5$, as
	$N\to\infty$,
	\begin{equation}\label{eq:Var-TVCc}
		\Var(d_{\mathbf p}^{c})
		=\Bigl(1-\tfrac2\pi\Bigr)\sum_{i=1}^{n}\sigma_i^2
		+\tfrac4\pi\!\!\sum_{1\le i<k\le n}\!\!\sigma_i\sigma_k\,\Phi(\rho_{i,k})+O(1),
	\end{equation}
	where $\Phi$ is the Nabeya folded-normal covariance factor (Sheppard's arcsine
	formula~\cite{sheppard1899}, in this form~\cite{nabeya1951}) used in
	\cite[Theorem~3.3]{banic_elezovic_p4}, and $\rho_{i,k}$ is defined by
	\eqref{eq:rho-cyclic}.
\end{proposition}

	The cyclic pairs consist of $n$ adjacent pairs, including the wrap-junction pairs
	$(n,1)$ and $(n-1,n)$, and $n(n-3)/2$ disjoint pairs. The $O(1)$ remainder is
	understood at fixed $n$ and fixed interior $\mathbf p$; its implicit constant is a
	sum over $\asymp n^2$ pairs and is not uniform in $n$.

\begin{proof}
	Each $|\tilde D_i|$ and each pair covariance is a property of a small marginal that
	is independent of whether the chain is open or closed (the transfer principle,
	Remark~\ref{rem:cyc-transfer}): the per-edge term is
	$\Var(|\tilde D_i|)=(1-\tfrac2\pi)\sigma_i^2+O(1)$ by the per-edge calculation
	in \cite[Theorem~3.3]{banic_elezovic_p4}. Each matched pair is
	\(\Cov(|\tilde D_i|,|\tilde D_k|)
	=\tfrac2\pi\sigma_i\sigma_k\Phi(\rho_{i,k})+O(1)\)
	by \cite[Proposition~3.6]{banic_elezovic_p4}; the assumption \(n\ge5\) leaves positive mass
	outside every three- or four-bin pair, exactly as that proposition requires. The
	wrap edge $D_n=X_1-X_n$ contributes its own per-edge term, two adjacent covariances
	(its junction triples $(X_{n-1},X_n,X_1)$ and $(X_n,X_1,X_2)$), and disjoint
	covariances with the remaining edges; summing over the cyclic multiplicities
	(Lemma~\ref{lem:cyc-mult}) gives \eqref{eq:Var-TVCc}.
\end{proof}

\begin{remark}[The deferred $O(1)$ remainder]\label{rem:TVCc-deferred}
	The leading sum is elementary --- every term a Nabeya value. What remains
	open is the $O(1)$ remainder, whose per-edge and adjacent constants are the
	\emph{non-uniform} cyclic analogues $c_L(\mathbf q),a_*(\mathbf q)$ of the uniform
	$c_L,a_*$ --- the same constants deferred for the linear $d_{\mathbf p}$ in
	Part~4. The general matched covariance is now in closed form
	\cite[Proposition~5.3]{banic_elezovic_p4}, but the \(O(1)\) correction \emph{beyond} its
	Gaussian value, like the uniform $a_*$, is not; the wrap edge adds nothing new to that
	difficulty, only its own term and two adjacent folded covariances.
\end{remark}

\section{Conclusion}
	\label{ssec:conclusion}

	This paper isolates the effect of closing the open chain of multinomial bins into a
	cycle. Under the uniform null the leading cyclic variance is independent of $n$ and
	equals the limiting open-chain value $\sigma_\infty^2$, while the first correction is
	the exact linear term $\beta_c(n)=\beta_1 n$. The proof uses no new marginal
	asymptotics: the cyclic statistic recombines the same per-edge, adjacent, and
	disjoint blocks as the open-chain statistic, but with cyclic multiplicities. The
	boundary cases $n=3$ and $n=4$ fit this beta-level formula without exception; at
	the next order only the deterministic four-bin total at $n=4$ leaves a parity
	oscillation.

	At second order the same recombination gives the conditional closed form
	\[
		\gamma_c(n)=\frac{(64a-3)n^2-4n}{32\pi}\qquad(n\ge5),
	\]
	where the per-edge and disjoint contributions are proved. The adjacent constant is
	$a=\tfrac1{32}-\tfrac{4\sqrt3}{243}$; it is identified to high numerical precision
	and is consistent with the assembled cyclic variance. A complete first-principles
	proof of this lattice coefficient remains the only conditional point in the
	formula. Uniform remainder control is also the natural missing ingredient for
	turning the asymptotic no-mode statement into a global fixed-$N$ result.

	The main consequence of the transfer identity is for the \emph{linear} statistic:
	comparison with the mode-free symmetric cyclic profile explains the traveling mode of
	the linear variance, and uniform control of the three ridge blocks proves its local
	optimizer $K_0=C_1/\beta_1=16.0429\ldots$ on every compact ridge window, upgrading the
	previously heuristic $\sqrt N$ mode law of the comb test to a theorem. The one global
	step left open is that the exact integer maximizer cannot escape this window.

	For fixed alternatives the cyclic statistic behaves exactly as expected from the
	added wrap edge. In the fully mismatched reference-centered statistic the cyclic variance is
	an explicit quadratic form in the edge means. For the matched one-sample statistic
	$d_{\mathbf p}^{c}$, the leading non-uniform variance is the cyclic
	Nabeya sum of Proposition~\ref{prop:Var-TVCc}. Its $O(1)$ refinement is
	governed by the same non-uniform constants as in the open-chain variance theory of
	Part~4~\cite{banic_elezovic_p4}; the passage from an open chain to a cycle changes
	the multiplicities and adds the wrap edge, but does not introduce a new local
	analytic obstruction.

\appendix

\section{Uniform control of the cyclic ridge blocks}
	\label{app:ridge-blocks}

	This appendix supplies the two uniform estimates used in
	\S\ref{ssec:uniform-blocks}: the aperiodicity and uniform spectral gap in the
	disjoint Fourier block, and the averaging of the deterministic adjacent
	asymptotic over a binomial three-bin total. Throughout this appendix,
	``uniformly'' means that the implicit constants may depend on the fixed compact
	ridge window, equivalently on a fixed \(C_0\) in \(n\le C_0\sqrt N\), but not on
	\(n,N\).

\begin{lemma}[Binomial mixing tails]\label{lem:app-mixing-tails}
	Let \(M\sim\Bin(N,d/n)\), \(2\le d\le4\), and let
	\(\bar M_d=dN/n\). If \(3\le n\le C_0\sqrt N\), then, for each fixed \(r>0\),
	\[
		\Pr\bigl(M\le\tfrac12\bar M_d\bigr)\le e^{-c_r\bar M_d},\qquad
		\E\bigl[M^{-1};M\ge\tfrac12\bar M_d\bigr]\le \frac{2}{\bar M_d},
	\]
	and
	\[
		\E\bigl[M^{-1/2};M\ge1\bigr]=O(\bar M_d^{-1/2}).
	\]
\end{lemma}

\begin{proof}
	The first estimate is Chernoff's binomial lower-tail bound; since
	\(\bar M_d=dN/n\ge(d/C_0)\sqrt N\), it is uniform in the ridge range. The
	second estimate is immediate from \(M^{-1}\le2/\bar M_d\) on the indicated event.
	For the third, split at \(M=\bar M_d/2\): on the bulk
	\(M^{-1/2}\le(2/\bar M_d)^{1/2}\), and the lower tail is exponentially small.
\end{proof}

\subsection*{The disjoint block}

	The exact Fourier reduction used in \cite[Theorem~4.6]{banic_elezovic_p2}
	contains, for the first Fourier block of the disjoint skeleton, the integral
	\begin{equation}\label{eq:app-JN1}
		J_{N,1}=\frac1{(2\pi)^2}\iint_{[-\pi,\pi]^2}
		\varphi(\theta)^{N-1}\,d\theta,
		\qquad
		\varphi(\theta)=1-\frac2n\bigl[(1-\cos\theta_1)+(1-\cos\theta_3)\bigr].
	\end{equation}
	The only point requiring care uniformly in the ridge is that this is a lattice
	characteristic function.

\begin{lemma}[Aperiodicity and uniform spectral gap]\label{lem:app-gap}
	Let
	\[
		S(\theta)=(1-\cos\theta_1)+(1-\cos\theta_3),\qquad
		\varphi(\theta)=1-\frac2n S(\theta).
	\]
	For \(n\ge5\), \(|\varphi(\theta)|=1\) on \([-\pi,\pi]^2\) if and only if
	\(\theta=0\). Moreover, for every \(\rho\in(0,\pi]\) there exists
	\(c=c(\rho)>0\), independent of \(n,N\), such that
	\[
		|\varphi(\theta)|\le1-\frac cn
		\qquad(\rho\le\|\theta\|_\infty\le\pi).
	\]
	Consequently, uniformly for \(5\le n\le C_0\sqrt N\),
	\[
		\int_{\rho\le\|\theta\|_\infty\le\pi}|\varphi(\theta)|^{N-1}\,d\theta
		=O(e^{-c'\sqrt N})
	\]
	for some \(c'>0\).
\end{lemma}

\begin{proof}
	Since \(0\le S\le4\), one has \(\varphi\le1\). Equality \(\varphi=1\) forces
	\(S=0\), hence \(\theta=0\). Equality \(\varphi=-1\) would force \(S=n\), which is
	impossible for \(n\ge5\). Thus the origin is the only unimodular point. For the
	gap, set \(K_\rho=\{\rho\le\|\theta\|_\infty\le\pi\}\). The continuous function
	\(S\) is positive on \(K_\rho\), so \(s_\rho=\min_{K_\rho}S>0\). If
	\(\varphi\ge0\), then
	\(1-|\varphi|=1-\varphi=2S/n\ge2s_\rho/n\). If \(\varphi<0\), then
	\[
		1-|\varphi|=1-\Bigl(\frac2nS-1\Bigr)=2-\frac2nS
		\ge2-\frac8n\ge\frac2n,
	\]
	because \(n\ge5\). This proves the gap with \(c=\min(2s_\rho,2)\). The exponential
	tail bound follows from \((1-c/n)^{N-1}\le e^{-c(N-1)/n}\) and
	\(n\le C_0\sqrt N\).
\end{proof}

\begin{lemma}[Uniform disjoint bound]\label{lem:app-disj}
	Uniformly for \(5\le n\le C_0\sqrt N\),
	\[
		\Cov(|D_1|,|D_3|)=-\frac1{\pi n}+O(1/N).
	\]
\end{lemma}

\begin{proof}
	By Lemma~\ref{lem:app-gap}, the part of \eqref{eq:app-JN1} outside
	\(\|\theta\|_\infty<\rho\) is exponentially small. On this neighborhood set
	\(\theta=\sqrt{n/(N-1)}\,u\). Since
	\[
		\varphi(\theta)=1-\frac{\theta_1^2+\theta_3^2}{n}
		+O\!\left(\frac{\theta_1^4+\theta_3^4}{n}\right),
	\]
	we get
	\[
		(N-1)\log\varphi(\theta)
		=-(u_1^2+u_3^2)+O\!\left(\frac nN |u|^4\right),
	\]
	uniformly for \(n\le C_0\sqrt N\). The Jacobian is \(n/(N-1)\), and dominated
	Laplace expansion gives
	\[
		J_{N,1}
		=\frac1{(2\pi)^2}\frac n{N-1}\iint_{\mathbb R^2}e^{-(u_1^2+u_3^2)}\,du
		\left[1+O\!\left(\frac nN\right)\right]
		=\frac{n}{4\pi(N-1)}
		\left[1+O\!\left(\frac nN\right)\right].
	\]
	Thus the \(k=1\) term in \cite[Theorem~4.6]{banic_elezovic_p2} is
	\[
		-\frac{4N}{n^2}J_{N,1}
		=-\frac1{\pi n}+O(1/N).
	\]
	The remaining \(k\ge2\) block in the same theorem is bounded by a constant times
	\[
		\frac{N^2}{n^2}
		\left[
		\int_{-\pi}^{\pi} a(\theta)e^{-\gamma(N-2)a(\theta)}\,d\theta
		\right]^2,
	\]
	with \(a(\theta)\asymp\theta^2/n\). The one-dimensional integral is
	\(O(n^{1/2}N^{-3/2})\), so this block is \(O((nN)^{-1})=O(N^{-1})\).
\end{proof}

\subsection*{The adjacent block}

	The adjacent block uses the deterministic adjacent skeleton from
	\cite[Theorem~5.3]{banic_elezovic_p2}:
	\begin{equation}\label{eq:app-hadj}
		\hadj(\ell)=c_L\ell+a_*+O(\ell^{-1}),\qquad
		c_L=\frac{2\sqrt3}{3\pi}+\frac19,\quad
		a_*=-\frac{2\sqrt3}{9\pi}.
	\end{equation}
	The \(O(\ell^{-1})\) remainder is a deterministic bound for a one-variable
	sequence, hence is uniform in \(\ell\).

\begin{lemma}[Averaging the adjacent skeleton]\label{lem:app-adj}
	Let \(M=X_1+X_2+X_3\sim\Bin(N,3/n)\). Uniformly for
	\(3\le n\le C_0\sqrt N\),
	\[
		\E\,\hadj(M)=c_L\,\frac{3N}{n}+a_*+O(n/N).
	\]
	Consequently,
	\[
		\Cov(|D_1|,|D_2|)
		=\frac{2B}{n}N+a_*+\frac{n-2}{2\pi n}+O(n/N).
	\]
\end{lemma}

\begin{proof}
	Put \(\bar M_3=3N/n\) and \(\ell_0=\bar M_3/2\). On \(M\ge\ell_0\), the
	remainder in \eqref{eq:app-hadj} is \(O(M^{-1})\le O(\ell_0^{-1})\), so by
	Lemma~\ref{lem:app-mixing-tails},
	\[
		\E[\hadj(M);M\ge\ell_0]
		=c_L\E[M;M\ge\ell_0]+a_*\Pr(M\ge\ell_0)+O(n/N).
	\]
	On the complementary event we use the trivial bound \(0\le\hadj(m)\le m^2\le N^2\)
	and the exponentially small lower tail. Hence
	\[
		\E\,\hadj(M)=c_L\,\frac{3N}{n}+a_*+O(n/N).
	\]
	Finally, Lemma~\ref{lem:peredge-uniform} gives
	\[
		(\E|D_1|)^2=\frac{4N}{\pi n}-\frac{n-2}{2\pi n}+O(n/N),
	\]
	and subtracting this from \(\E\,\hadj(M)\), using
	\(3c_L-4/\pi=2B\), proves the covariance formula.
\end{proof}

\begin{remark}[What the appendix closes]\label{rem:app-closed}
	Lemma~\ref{lem:app-gap} closes the aperiodicity and spectral-gap step in the
	disjoint block, and Lemma~\ref{lem:app-adj} justifies the binomial averaging of
	the adjacent asymptotic. Thus the ridge-local expansion in
	Theorem~\ref{thm:K0-ridge} is unconditional at the level needed for \(K_0\).
	The global localization of the exact integer maximizer and the first-principles
	proof of the second-order triangular-lattice constant \(a\) remain separate
	questions.
\end{remark}

\section*{Acknowledgments}
	The authors used an AI assistant (Anthropic's Claude) for literature search, symbolic
	and numerical verification, and editorial support; all mathematical results were
	checked independently by the authors. The authors declare no competing interests.

\end{document}